\documentclass{gtart_h}  

\def\ifplaintex{\expandafter\ifx\csname documentclass\endcsname\relax}

\def\ifplaintex{\expandafter\ifx\csname documentclass\endcsname\relax}

\ifplaintex 
\else
\fi

\expandafter\ifx\csname epsfbox\endcsname\relax\input epsf\fi

\def\gt{{\mathsurround=0pt\it $\cal G\mskip-2mu$eometry \&\ 
$\cal T\!\!$opology}}        

\def\gtp{{\mathsurround=0pt\it $\cal G\mskip-2mu$eometry \&\ 
$\cal T\!\!$opology $\cal P\!$ublications}}  

\def\lognumber#1{\def\thelognumber{#1}}
\def\volumenumber#1{\def\thevolumenumber{#1}}
\def\papernumber#1{\def\thepapernumber{#1}}
\def\volumeyear#1{\def\thevolumeyear{#1}}

\def\pagenumbers#1#2{\def\startpage{#1}\def\finishpage{#2}}
\def\published#1{\def\publishdate{#1}}
\def\proposed#1{\def\theproposer{#1}}
\def\seconded#1{\def\theseconders{#1}}
\def\received#1{\def\receiveddate{#1}}
\def\revised#1{\def\reviseddate{#1}}
\def\accepted#1{\def\accepteddate{#1}}
\def\asciititle#1{\def\theasciititle{#1}}

\def\coverauthors#1{\def\thecoverauthors{#1}}
\def\asciiauthors#1{\def\theasciiauthors{#1}}

\def\asciiemail#1{\def\theasciiemail{#1}}

\long\def\asciiabstract#1{\long\def\theasciiabstract{#1}}
\def\asciikeywords#1{\def\theasciikeywords{#1}}

\let\\\par\let\thelognumber\relax
\let\thevolumenumber\relax\let\thepapernumber\relax
\let\thevolumeyear\relax\let\thesamplenumber\relax\let\startpage\relax
\let\finishpage\relax\let\publishdate\relax\let\receiveddate\relax
\let\reviseddate\relax\let\accepteddate\relax\let\theasciititle\relax
\let\theasciiauthors\relax
\let\theasciiabstract\relax\let\theasciikeywords\relax
\let\theasciiemail\relax\let\theshortauthors\relax\let\theshorttitle\relax
\let\thecoverauthors\relax

\long\def\maketitlep{   

\count0=\startpage

\gt\hfill      
\hbox to 77pt{\vbox to 0pt{\vglue -15pt\epsfbox{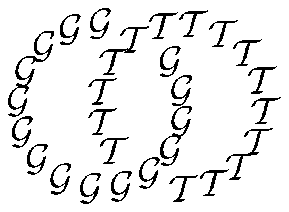}\vss}\hss}
\break
{\small\ifx\thesamplenumber\relax 
Volume \else Sample
\fi\thevolumenumber\ (\thevolumeyear)
\startpage--\finishpage\nl
Published: \publishdate}
\vglue 0.5truein plus 0.4fil minus 0.1truein

{\parskip=0pt\leftskip 0pt plus 1fil\def\\{\par\smallskip}{\ifplaintex\large
\else\Large\fi\bf\thetitle}\par\medskip}   

\vglue 0pt plus 0.1fil 

{\parskip=0pt\leftskip 0pt plus 1fil\def\\{\par}{\sc\theauthors}
\par\medskip}

\vglue 0pt plus 0.1fil 

{\small\parskip=0pt\let\newline\\
{\leftskip 0pt plus 1fil\def\\{\par}{\sl\theaddress}\par}
\expandafter\ifx\theemail\relax    
\relax\else\vglue 5pt plus 0.02fil minus 2pt\def\\{\stdspace{\rm 
and}\stdspace} 
\cl{Email:\stdspace\tt\theemail}\fi
\ifx\theurl\relax                  
\relax\else\vglue 5pt plus 0.02fil minus 2pt\def\\{\stdspace{\rm 
and}\stdspace}
\cl{URL:\stdspace\tt\theurl}\fi\par}

\vglue 7pt plus 0.3fil minus 3pt

{\bf Abstract}
\vglue 5pt plus 0.1fil minus 2pt

\theabstract

\vglue 7pt plus 0.3fil minus 3pt

{\bf AMS Classification numbers}\quad Primary:\quad \theprimaryclass

Secondary:\quad \thesecondaryclass

\vglue 5pt plus 0.3fil minus 2pt

{\bf Keywords:}\quad \thekeywords

\vglue 10pt plus 0.5fil minus 5pt

{\small  Proposed: \theproposer\hfill Received: \receiveddate\nl
Seconded: \theseconders\hfill 
\ifx\reviseddate\relax                         
Accepted: \accepteddate                        
\else
Revised: \reviseddate                          
\fi}
\eject
}       

\font\phead=cmsl9 scaled 950
\font=cmsl9 scaled 1050
\font\pnum=cmbx10 scaled 913
\font\lnum=cmbx10 
\font\pfoot=cmsl9 scaled 950
\font=cmsl9 scaled 1050
\ifplaintex
\headline{\vbox to 0pt{\vskip -4.5mm\line{\small\phead\ifnum
\count0=\startpage ISSN 1364-0380 (on line)
1465-3060 (printed) \hfill {\pnum\folio}\else\ifodd\count0\def\\{ }%
\ifx\theshorttitle\relax\thetitle\else\theshorttitle\fi\hfill{\pnum\folio}
\else\def\\{ and }{\pnum\folio}\hfill\ifx\theshortauthors\relax\theauthors
\else\theshortauthors\fi\fi\fi}\vss}}
\footline{\vbox to 0pt{\vglue 0mm\line{\small\pfoot\ifnum\count0=\startpage
\copyright\ \gtp\hfill\else
\gt, Volume \thevolumenumber\ (\thevolumeyear)\hfill\fi}\vss
}}
\else
\makeatletter
\def\@oddhead{{\small\lhead\ifnum\count0=\startpage ISSN 1364-0380 (on line)
1465-3060 (printed) \hfill {\lnum\number\count0}\else\ifodd\count0
\def\\{ }\ifx\theshorttitle\relax \thetitle \else\theshorttitle\fi\hfill
{\lnum\number\count0}\else\def\\{ and }{\lnum\number\count0}
\hfill\ifx\theshortauthors\relax 
\theauthors\else\theshortauthors\fi\fi\fi}}\def\@evenhead{\@oddhead}
\def\@oddfoot{\small\lfoot\ifnum\count0=\startpage\copyright\ \gtp\hfill\else
\gt, Volume \thevolumenumber\ (\thevolumeyear)\hfill\fi}
\def\@evenfoot{\@oddfoot}
\makeatother
\fi

\newwrite\gtoutfile
\long\gdef\makeheadfile{  
{\def\\{, }\def\s{ }
\immediate\openout\gtoutfile head.xxx
\immediate\write\gtoutfile{Proxy-for: \ifx\theasciiauthors\relax
\theauthors\else\theasciiauthors\fi\s<\ifx\theasciiemail\relax\theemail\else\theasciiemail\fi>}
\immediate\write\gtoutfile{\noexpand\\}
\immediate\write\gtoutfile{Authors: \ifx\theasciiauthors\relax
\theauthors\else\theasciiauthors\fi}
{\def\\{ }\immediate\write\gtoutfile{Title: \ifx\theasciititle\relax
\thetitle\else\theasciititle\fi}}
\immediate\write\gtoutfile{Subj-class: GT or SG or MG etc}
\immediate\write\gtoutfile{MSC-class: \theprimaryclass\ifx\thesecondaryclass\relax\else, \thesecondaryclass\fi}
\immediate\write\gtoutfile{Journal-ref: Geom. Topol. \thevolumenumber
(\thevolumeyear) \startpage-\finishpage}
\immediate\write\gtoutfile{Comments: Published by Geometry and Topology at}
\immediate\write\gtoutfile{\s\s http://www.maths.warwick.ac.uk/gt/GTVol\thevolumenumber/paper\thepapernumber.abs.html}
\immediate\write\gtoutfile{\noexpand\\}
\immediate\write\gtoutfile{}
\ifx\theasciiabstract\relax
\immediate\write\gtoutfile{\theabstract}\else
\immediate\write\gtoutfile{\theasciiabstract}\fi
\immediate\write\gtoutfile{}
\immediate\write\gtoutfile{\noexpand\\}
\immediate\write\gtoutfile{}
\immediate\closeout\gtoutfile}}  

\def\maketitlepage{\maketitlep\makeheadfile}
\let\maketitle\maketitlepage

\lognumber{509}
\received{28 October 2004}
\volumenumber{9}\papernumber{29}\volumeyear{2005}
\pagenumbers{1253}{1293}   
\revised{20 July 2005}
\published{24 July 2005}
\accepted{3 July 2005}
\proposed{Walter Neumann}
\seconded{Joan Birman, Vaughan Jones}

\usepackage{amssymb,epic,eepic,epsfig,amsmath}

\theoremstyle{plain}

\newtheorem{theorem}{Theorem}
\newtheorem{proposition}{Proposition}[section]
\newtheorem{lemma}[proposition]{Lemma}
\newtheorem{corollary}[proposition]{Corollary}

\theoremstyle{definition}
\newtheorem{definition}[proposition]{Definition}

\theoremstyle{remark}

\newtheorem{remark}[proposition]{Remark}

\newcommand{\psdraw}[2]
         {\begin{array}{c} \hspace{-1.3mm}
    \raisebox{-4pt}{\epsfig{figure=draws/#1.eps,width=#2}}
    \hspace{-1.9mm}\end{array}}

\newlength{\standardunitlength}
\def\KK{{\mathcal K}}
\def\BN{\mathbb N}
\def\BZ{\mathbb Z}

\def\BQ{\mathbb Q}
\def\BR{\mathbb R}
\def\BC{\mathbb C}

\def\BK{\mathbb K}
\def\A{\mathcal A}
\def\B{\mathcal B}

\def\DD{\mathcal D}
\def\BB{\mathcal B}
\def\K{\mathbf K}

\def\T{\mathcal T}

\def\N{{\mathbb N}}

\def\T{\mathcal T}
\def\longto{\longrightarrow}
\def\msl{\mathfrak{sl}}
\def\la{\langle}
\def\ra{\rangle}

\def\calI{\mathcal I}

\def\fg{\mathfrak{g}}

\newcommand{\tr}{\operatorname{tr}}

\newcommand{\id}{\operatorname{id}}
\newcommand{\qbinom}[2]{\text{$\left[\begin{array}{c}#1\\ #2\end{array}
\right]$}}

\def\a{\alpha}

\def\l{\lambda}

\def\s{\sigma}

\def\ga{\gamma}

\def\j{{\mathbf j}}
\def\m{{\mathbf m}}

\def\e{\mathbf e}

\def\d{\delta}
\def\b{\beta}

\def\s{\sigma}

\def\al{\alpha}

\def\Z{\BZ}

\def\hh{\mathfrak{h}}

\def\bk{\mathbf{k}}
\def\Weyl{W}
\def\UU{{\mathcal U}}
\def\Q{\BQ}
\def\n{{\mathbf n}}
\def\f{{\mathbf f}}
\def\ep{\epsilon}
\def\tr{\mathrm{tr}}
\def\bs#1{\boldsymbol{#1}}

\begin{document}


\title{The colored Jones function is $q$--holonomic}
\asciititle{The colored Jones function is q-holonomic}

\author{Stavros Garoufalidis\\Thang T\,Q L\^e}
\coverauthors{Stavros Garoufalidis\\Thang TQ L\noexpand\^e}
\asciiauthors{Stavros Garoufalidis\\Thang TQ Le}
\address{School of Mathematics, Georgia Institute of Technology\\
         Atlanta, GA 30332-0160, USA}

\gtemail{\mailto{stavros@math.gatech.edu}, \mailto{letu@math.gatech.edu}}
\asciiemail{stavros@math.gatech.edu, letu@math.gatech.edu}
\urladdr{http://www.math.gatech.edu/~stavros} 

\asciikeywords{Holonomic functions, Jones polynomial, Knots,
WZ algorithm, quantum invariants, D-modules, multisums, hypergeometric
functions}

\begin{abstract}
A function of several variables is called holonomic if, roughly
speaking, it is determined from finitely many of its values via finitely
many linear recursion relations with polynomial coefficients. 
Zeilberger was the first to notice that
the abstract notion of holonomicity can be applied to verify, in a
systematic and computerized way, combinatorial identities among
special functions. Using a general state sum definition of the
colored Jones function of a link in 3--space, we prove from first
principles that the colored Jones function is a multisum of
a $q$--proper-hypergeometric function, and thus it is $q$--holonomic.
We demonstrate our results by computer calculations.
\end{abstract}
\asciiabstract{%
A function of several variables is called holonomic if, roughly
speaking, it is determined from finitely many of its values via
finitely many linear recursion relations with polynomial coefficients.
Zeilberger was the first to notice that the abstract notion of
holonomicity can be applied to verify, in a systematic and
computerized way, combinatorial identities among special
functions. Using a general state sum definition of the colored Jones
function of a link in 3-space, we prove from first principles that
the colored Jones function is a multisum of a
q-proper-hypergeometric function, and thus it is q-holonomic.
We demonstrate our results by computer calculations.}

\maketitle


\section{Introduction}
\label{sec.intro}

\subsection{Zeilberger meets Jones}
\label{sub.ZmeetsJ}

The colored Jones function of a framed knot $\KK$ in 3--space
$$
J_\KK\co \BN\longto\BZ[q^{\pm 1/4}]
$$ is a sequence of Laurent
polynomials that essentially measures the Jones polynomial of a
knot and its cables. This is a powerful but not well understood
invariant of knots. As an example, the colored Jones function of
the $0$--framed right-hand trefoil is given by
$$
J_\KK(n) = \frac{q^{1/2-n/2}}{1-q^{-1}} \sum_ {k=0}^{n-1} q^{-k n}
(1- q^{-n})(1-q^{1-n}) \dots (1 - q^{k-n}).
$$
Here $J_\KK(n)$ denotes the Jones polynomial of the $0$--framed knot
$\KK$ colored by the $n$--dimensional irreducible representation of $\msl_2$, 
and normalized by 
$J_{\text{unknot}}(n)=(q^{n/2}-q^{-n/2})/(q^{1/2}-q^{-1/2})$.

Only a handful of 
knots have such a simple formula. However, as we shall see all knots
have a {\em multisum}
 formula. Another way to look at the colored Jones function
of the trefoil is via the following 3--term recursion formula:
$$
J_\KK(n)= \frac{q^{n-1} + q^{4-4n}- q^{-n} -
q^{1-2n}}{q^{1/2}(q^{n-1}-q^{2-n} )} J_\KK(n-1) + \frac{q^{4-4n}
-q^{3-2n}}{q^{2-n} -q^{n-1}}J_\KK(n-2)
$$
with initial conditions: $J_\KK(0)=0, \quad J_\KK(1)=1$.

In this paper we prove that the colored Jones function of any knot
satisfies a linear recursion relation, similar to the above
one. For a few knots this was obtained by Gelca and his colleagues
\cite{Gelca1,Gelca2}. (In \cite{Gelca1} a more complicated  5--term
recursion formula for the trefoil was established).

Discrete functions that satisfy a nontrivial difference recursion
relation are known by another name: they are $q$--holonomic.

Holonomic functions were introduced by IN Bernstein \cite{B1,B2}
and M Saito. The latter coined the term holonomic, that is
a function which is entirely determined by the law
of its differential equation, together with finitely many
initial conditions. Bernstein used holonomic functions to prove a
conjecture of Gelfand on the analytic continuation of operators.
Holonomicity and the related notion of $D$--modules are a tool in
studying linear differential equations from the point of view of
algebra (differential Galois theory), algebraic geometry, and
category theory. For an excellent introduction on holonomic
functions and their properties, see \cite{Bo} and \cite{Cou}.

Our approach to the colored Jones function owes greatly to
Zeilberger's work. Zeilberger noticed that
the abstract notion of holonomicity can be applied to verify, in a systematic
and computerized way, combinatorial identities among special functions,
\cite{Z} and also \cite{WZ, PWZ}.

A starting point for Zeilberger, the so-called {\em operator approach},
is to replace functions by the recursion relations that they satisfy.
This idea leads in a natural way
to noncommutative algebras of operators that act on a function, together with
left ideals of annihilating operators.

To explain this idea concretely, consider the operators $E$ and
$Q$ which act on a {\em discrete function} (that is, a function of a discrete
variable $n$) $f\co \BN\longto\BZ[q^{\pm}]$ by:
$$
(Q f)(n)=q^{n}f(n) \hspace{1cm}
(E f)(n) = f(n+1).
$$
It is easy to see that $EQ=qQE$, and that $E,Q$ generate a noncommutative
$q$--{\em Weyl algebra} generated by noncommutative polynomials in $E$
and $Q$, modulo the relation $EQ=qQE$:
$$
\A=\BZ[q^{\pm}]\la Q,E\ra/(EQ=qQE)
$$
Given a discrete function $f$ as above, consider the {\em recursion ideal}
$\calI_f=\{P \in \A \, | Pf=0\}$. It is easy to see that it is a left
ideal of the $q$--Weyl algebra. We say that $f$ is $q$--{\em holonomic} iff
$\calI_f \neq 0$.

In this paper we prove that:

\begin{theorem}
\label{thm.1}
The colored Jones function of every knot is $q$--holonomic.
\end{theorem}

Theorem \ref{thm.1} and its companion Theorem \ref{thm.c}
are effective, as their proof reveals.

\begin{theorem}
\label{thm.c}$\phantom{99}$

{\rm (a)}\qua
The $E$--order of the colored Jones function of a knot is bounded above by
an exponential function in the number of crossings. 

{\rm (b)}\qua
For every knot $\KK$ there exist a natural number $n(\KK)$, such that
$n(\KK)$ initial values of the colored Jones function determine
the colored Jones function of $\KK$. In other words, the colored Jones function
is determined by a finite list. $n(\KK)$ is bounded above by an exponential
function in the number of crossings.
\end{theorem}

Computer calculations are given in Section \ref{sec.computer}.
In relation to (b) above, notice that the $q$--Weyl algebra is
{\em noetherian}; thus every left ideal is finitely generated. The theorem
states more, namely that the we can compute (via elimination) a basis
for the recursion ideal of the colored Jones function of a knot.

Let us end the introduction with some remarks.

\begin{remark}
\label{rem.i2} The colored Jones function can be defined for every
simple Lie algebra $\fg$. Our proof of Theorem \ref{thm.1} generalizes
and proves that the $\fg$--colored Jones function of a knot is
$q$--holonomic (except for $G_2$), see Theorem \ref{mainn} below.
\end{remark}

\begin{remark}
\label{rem.ilink} The colored Jones function can be defined for colored
links in 3--space. Our proof of Theorem \ref{thm.1} proves that the
colored Jones function of a link is $q$--holonomic in all variables,
see Section \ref{sub.thm.1}.
\end{remark}

\begin{remark}
\label{rem.i3}
It is well known that computing $J(n)$ for any fixed $n>1$
is a \#{\em P--complete} problem. Theorem \ref{thm.1} claims that this
sequence of \#P--complete problems is no worse than any of its terms.
\end{remark}

\begin{remark}
\label{rem.i5}
The proof of Theorem \ref{thm.1} indicates that many statistical mechanics
models, with complicated partition functions that depend on several variables,
are holonomic, provided that their local weights are holonomic. This
observation may be of interest to statistical mechanics.
\end{remark}

\subsection{Synonymous notions to holonomicity}
\label{sec.synonymous}

We have chosen to phrase the results of our paper mostly using the 
high-school language of linear recursion relations. We could have used
synonymous terms such as linear $q$--difference equations, or
$q$--holonomic functions, or $D$--modules, or maximally overdetermined systems
of linear PDEs which is more common in the area of algebraic analysis,
see for example \cite{Ml}. The geometric notion of $D$--modules gives 
rise to geometric invariants of knots, such as the characteristic variety
introduced by the first author in \cite{G2}. The characteristic variety
is determined by the colored Jones function of a knot and is conjectured to
be isomorphic to the $\mathfrak{sl}_2(\BC)$--character variety of a knot, 
viewed from the boundary torus.
This, so-called {\em AJ Conjecture}, formulated by the first author is known
to hold for all torus knots (due to Hikami, \cite{Hi}), 
and infinitely many 2--bridge knots (due to the second-author, \cite{Le2}).

Thus, there is nontrivial geometry encoded in the linear recursion relations
of the colored Jones function of a knot.

\subsection{Plan of the proof}
\label{sub.plan}

In Section \ref{sec.qholo}, we discuss in detail the notion of a $q$--holonomic
function. We give examples of $q$--holonomic functions (our building blocks),
together with rules that create $q$--holonomic functions from known
ones.

In Section  \ref{sec.state}, we discuss the colored Jones function
of a link in $3$--space, using state sums associated to a planar
projection of the link. The colored Jones function is built out of
local building blocks (namely, $R$--matrices) associated to the
crossings, which are assembled together in a way dictated by the
planar projection. The main observation is that the $R$--matrix is
$q$--holonomic in all variables, and that the assembly preserves
$q$--holonomicity. Theorem \ref{thm.1} follows. As a bonus, we
present the colored Jones function as a multisum of a $q$--proper
hypergeometric function.

In Section \ref{sec.cyclo} we show that the cyclotomic function of a knot
(a reparametrization of the colored Jones function, introduced by Habiro,
with good integrality properties) is $q$--holonomic, too. We achieve this
by studying explicitly a change of basis for representations of $\msl_2$.

In Section \ref{sec.effective} we give a theoretical review about complexity
and computability of recursion relations of $q$--holonomic functions,
following Zeilberger.
These ideas solve the problem of finding recursion relations
of $q$--holonomic functions which are given by multisums of $q$ proper
hypergeometric functions. It is a fortunate coincidence (?) that the
colored Jones function can be presented by such a multisum, thus we can
compute its recursion relations. Theorem \ref{thm.c} follows.

Section \ref{sec.computer} is a computer implementation of the previous
section, where we use {\tt Mathematica} packages developed by A Riese.

In Section \ref{sec.lie} we discuss the $\fg$--colored Jones function
of a knot, associated to a simple Lie algebra $\fg$. Our goal is to
prove that the $\fg$--colored Jones function is $q$--holonomic in all
variables (see Theorem \ref{mainn}). In analogy with the $\fg=\msl_2$ case,
we need to show that the local building block, the $R$--matrix, is
$q$--holonomic in all variables. This is a trip to the world of quantum groups,
which takes up the rest of the section, and ends with an appendix which
computes (by brute-force) structure constants of quantized enveloping Lie
algebras in the rank $2$ case.

\subsection{Acknowledgements}

The authors wish to express their gratitude to D Zeilberger who
introduced them to the wonderful world of holonomic functions and its
relation to discrete mathematics. Zeilberger is a philosophical
co-author of the present paper.

In addition, the authors wish to thank A Bernstein, G Lusztig, and N
Xie for help in understanding quantum group theory, G Masbaum and K
Habiro for help in a formula of the cyclotomic expansion of the
colored Jones polynomial, and D Bar-Natan and A Riese for computer
implementations.

The authors were supported in part by National Science Foundation.

\section{$q$--holonomic and $q$--hypergeometric functions}
\label{sec.qholo}

Theorem \ref{thm.1} follows from the fact that the colored Jones
function can be built from elementary blocks that are
$q$--holonomic, and the operations that patch the blocks together
to give the colored Jones function preserve $q$--holonomicity.

IN Bernstein defined the notion of holonomic functions $f\co \BR^r
\longto \BC$, \cite{B1,B2}.  For an excellent and complete
account, see Bjork \cite{Bj}. Zeilberger's brilliant idea was to
link the abstract notion of holonomicity to the concrete problem
of algorithmically proving combinatorial identities among
hypergeometric functions, see \cite{Z, WZ} and also \cite{PWZ}. This
opened an entirely new view on combinatorial identities.

Sabbah extended Bernstein's approach to holonomic functions and
defined the notion of a $q$--holonomic function, see \cite{S} and
also \cite{C}.

\subsection{$q$--holonomicity in many variables}
We briefly review here the definition of $q$--holonomicity. First
of all, we need an $r$--dimensional version of the $q$--Weyl
algebra. Consider the operators $E_i$ and $Q_j$ for $1 \leq i,j
\leq r$ which act on discrete functions $f\co  \BN^r \longto
\BZ[q^{\pm}]$ by:
\begin{eqnarray*}
(Q_i f)(n_1,\dots,n_r)&=&q^{n_i}f(n_1,\dots, n_r) \\
(E_i f)(n_1,\dots,n_r) &=& f(n_1,\dots, n_{i-1}, n_i+1, n_{i+1}, \dots, n_r).
\end{eqnarray*}
It is easy to see that the following relations hold:
\begin{equation*}
\tag{$\text{Rel}_{q}$}
\begin{aligned}
Q_i Q_j&= Q_j Q_i  & E_i E_j &= E_j E_i \\
Q_i E_j&= E_j Q_i \, \text{for} \, i \neq j & E_i Q_i &= q Q_i E_i
\end{aligned}
\end{equation*}
We define the $q$--{\em Weyl algebra} $\A_r$ to be a noncommutative algebra
with presentation
$$
\A_r= \frac{\BZ[q^{\pm1}]\la Q_1, \dots, Q_r , E_1, \dots, E_r
\ra}{ (\text{Rel}_q)}.
$$
Given a discrete function $f$ with domain $\BN^r$ or $\BZ^r$ and
target space a $\BZ[q^{\pm 1}]$--module, one can define the left
ideal $\calI_f$ in $\A_r$ by
$$  
\calI_f := \{P \in \A_r | Pf=0 \}    .
$$ If we want to determine a function $f$ by a finite list of
initial conditions, it does not suffice to ensure that $f$
satisfies one nontrivial recursion relation if $r\ge 2$. The key
notion that we need instead is $q$--holonomicity.

Intuitively, a discrete function $f\co \BN^r \longto \BZ[q^{\pm}]$ is
$q$--holonomic if it satisfies a {\em maximally overdetermined
system} of linear difference equations with polynomial
coefficients. The exact definition of holonomicity is through
homological dimension, as follows.

Suppose $M= \A_r/I$, where $I$ is a left $\A_r$--module. Let $F_m$
be the sub-space of $\A_r$ spanned by polynomials in $Q_i,E_i$ of
total degree $\le m$. Then the module $\A_r/I$ can be approximated
by the sequence $F_m/(F_m\cap I), m=1,2,...$. It turns out that,
for $m >>1$,  the dimension (over the fractional field $\BQ(q)$)
of $F_m/(F_m\cap I)$ is a polynomial in $m$ whose degree $d(M)$ is
called the {\em homological dimension} of $M$.

Bernstein's {\em famous inequality} (proved by Sabbah in the
$q$--case, \cite{S}) states that $d(M) \geq r$, if $M\neq 0$ and $M$ has
{\em no monomial torsions}, ie, any non-trivial element of $M$ cannot be
annihilated by a monomial in $Q_i,E_i$. Note that the  left
$\A_r$--module $M_f:=\A_r \cdot f \cong \A_r/\calI_f$ does not
have monomial torsion.

\begin{definition}
\label{def.qholo}  We say that a discrete function $f$ is
$q$--holonomic if $d(M_f)\le r$.
\end{definition}

Note that if $d(M_f)\le r$, then by Bernstein's inequality, either
$M_f=0$ or $d(M_f)=r$. The former can happen only if $f=0$.

Although we will not use in this paper, let us point out
an alternative cohomological definition of dimension for 
a finitely generated $\A_r$ module $M$.
Let us define 
$$
c(M) := \min \{ j\in \BN \,\, | \,\, \mathrm{Ext}^j_{\A_r}(M, \A_r)\neq 0\}.
$$
Then the {\em homological dimension} $d(M) := 2r- c(M)$ equals to the 
dimension $d(M)$ defined above.

Closely related to $\A_r$ is the $q$--{\em torus algebra} $\T_r$ with
presentation
$$
\T_r= \frac{\BZ[q^{\pm 1}]\la Q_1^{\pm 1}, \dots, Q_r^{\pm 1} ,
E_1^{\pm 1}, \dots, E_r^{\pm 1} \ra}{ (\text{Rel}_q)}.
$$
Elements of $\T_r$ acts on the set of functions with domain
$\BZ^r$, but not on the set of functions with domain $\BN^r$. Note
that $\T_r$ is simple, but $\A_r$ is not. If $I$ is a left ideal
of $\T_r$ then the dimension of $\T_r/I$ is equal to that of
$\A_r/(I\cap A_r)$.

\subsection{Assembling $q$--holonomic functions}
\label{sub.assemble}

Despite the unwelcoming definition of $q$--holonomic functions, in this
paper we will use not the definition itself, but rather the {\em closure
properties} of the set of $q$--holonomic functions under some 
natural operations. 

{\bf Fact 0}

\begin{itemize}
\item
Sums and products of $q$--holonomic functions are $q$--holonomic.
\item
Specializations and extensions of $q$--holonomic functions are
$q$--holonomic. In other words, if $f(n_1,\dots,n_m)$ is
$q$--holonomic, the  so are the functions
\begin{align*} g(n_2,\dots,n_m):=&
f(a,n_2,\dots,n_m)\\\hbox{and}\qquad h(n_1,\dots,n_m,n_{m+1}):=&
f(n_1,\dots,n_m).\end{align*}
\item
Diagonals of $q$--holonomic functions are $q$--holonomic. In
other words, if $f(n_1,\dots,n_m)$ is $q$--holonomic, then  so is
the function
$$g(n_2,\dots,n_m) := f(n_2,n_2,n_3,\dots,n_m).$$
\item
Linear substitution. If $f(n_1,\dots,n_m)$ is $q$--holonomic, then  so is
the function, $g(n_1',\dots,n'_{m'})$, where each $n'_j$ is a
linear function of $n_i$.
\item
Multisums of $q$--holonomic functions are
$q$--holonomic. In other words, if $f(n_1,\dots,n_m)$ is
$q$--holonomic, the  so are the functions $g$ and $h$, defined by
\begin{align*}
g(a,b,n_2,\dots,n_m):=& \sum_{n_1=a}^b f(n_1,n_2,\dots,n_m)\\
h(a,n_2,\dots,n_m) :=& \sum_{n_1=a}^\infty f(n_1,n_2,\dots,n_m)
\end{align*}
(assuming that the latter sum is finite for each $a$).
\end{itemize}

For a user-friendly explanation of these facts and for many examples,
see \cite{Z,WZ} and \cite{PWZ}.

\subsection{Examples of $q$--holonomic functions}
Here are a few examples of $q$--holonomic functions. In fact, we
will encounter only sums, products, extensions, specializations,
diagonals, and multisums of these functions. In what follows we
usually extend the ground ring $\BZ[q^{\pm 1}]$ to the {\em fractional
field} $\BQ(q^{1/D})$, where $D$ is a positive integer. We also
use $v$ to denote a root of $q$, $v^2=q$.

For $n,k\in \BZ$, let
$$
\{n\}:= v^n -v^{-n},
\qquad [n]:= \frac{\{n\}}{\{1\}},
\qquad [n]! :=\prod_{i=1}^n
[i], \qquad \{n\}! := \prod_{i=1}^{n} \{i\}
$$
$$
\{n\}_k := \begin{cases}
\prod_{i=1}^k \{n-i+1\},  & \text{if $k\ge 0$}\\
0 & \text{if $k<0$} \end{cases}
$$
$$
\qbinom{n}{k} := \begin{cases}
\frac{\{n\}_k}{\{k\}_k}  & \text{if $k\ge 0$}\\
0 & \text{if $k<0$} \end{cases}.
$$
The first four functions are $q$--holonomic in $n$, and the last
two, as well as the delta function $\d_{n,k}$, are $q$--holonomic
in both $n$ and $k$.

\subsection{$q$--hypergeometric functions}
\label{sub.hyper}

\begin{definition}
\label{def.hyper}
A discrete function $f\co  \BZ^r \longto \BQ(q)$ is $q$--{\em hypergeometric}
iff $E_i f/f \in \BQ(q,q^{n_1}, \dots, q^{n_r})$ for all $i=1,\dots,r$.
\end{definition}

In that case, we know generators for the annihilation ideal of $f$.
Namely, let $E_if/f=(R_i/S_i)|_{Q_i=q^{n_i}}$
for $R_i,S_i \in \BZ[q,Q_1,\dots,Q_r]$.
Then, the annihilation ideal of $f$ is generated by $S_iE_i-R_i$.

All the functions in the previous subsections are
$q$--hypergeometric.

Unfortunately, $q$--hypergeometric functions are not always
$q$--holonomic. For example, $(n,k)\longto 1/[n^2+k^2]!$ is
$q$--hypergeometric but not $q$--holonomic. However,
$q$--proper-hypergeometric functions are $q$--holonomic.
The latter were defined by Wilf--Zeilberger as follows, \cite[Sec.3.1]{WZ}:

\begin{definition}
\label{def.proper}
A {\em proper $q$--hypergeometric} discrete
function is one of the form
\begin{equation}
\label{eq.proper}
F(n, \bk)=\frac{\prod_s (A_s;q)_{a_sn+\mathbf{b}_s. \bk+c_s}}{
\prod_t (B_t;q)_{u_t n+\mathbf{v}_t . \bk +w_t}}
q^{A(n,\bk)} 
\xi^{\bk}
\end{equation}
where $A_s,B_t \in \BK=\BQ(q)$, $a_s, u_t$ are integers,
$\mathbf{b}_s, \bk_s$ are vectors of $r$ integers, $A(n,\bk)$ is a
quadratic form, $c_s,w_s$ are variables and $\xi$ is an $r$ vector
of elements in $\BK$. Here, as usual
$$ (A;q)_n := \prod_{i=0}^{n-1} (1-Aq^i).$$
\end{definition}

\section{The colored Jones function for $\mathfrak{sl}_2$}
\label{sec.state}

\subsection{Proof of Theorem \ref{thm.1} for links}
\label{sub.thm.1}

We will formulate and prove  an analog of Theorem \ref{thm.1} (see
Theorem \ref{thm.link} below) for colored links.  Our proof will use
a {\em state-sum} definition of the colored Jones function, coming from
a representation of the quantum group $U_q(\msl_2)$, as was discovered
by Reshetikhin and Turaev in \cite{RT,Tu}.

Suppose $L$ is a framed, oriented link of $p$ components. Then the colored
Jones function $J_L\co  \BN^p \to \BZ[q^{\pm 1/4}] =\BZ[v^{\pm 1/2}]$ can
be defined  using  the representations of braid groups coming from
the quantum group $U_q(\msl_2)$.

\begin{theorem}
\label{thm.link}
The colored Jones function $J_L$ is $q$--holonomic.
\end{theorem}

\begin{proof}
We will present the definition of $J_L$ in the form most suitable
for us. Let $V(n)$ be the $n$--dimensional vector space over the
field $\BQ(v^{1/2})$ with basis
$\{e_0,e_1,\dots,e_{n-1}\}$, with $V(0)$ the zero vector space.

Fix a positive integer $m$. A linear operator
$$
A\co V(n_1)\otimes\dots \otimes V(n_m) \to V(n_1') \otimes\dots \otimes
V(n_m')
$$
can be described by the collection
$$
A_{a_1,\dots,a_m}^{b_1,\dots,b_m} \in \BQ(v^{1/2}),
$$
where
$$ A(e_{a_1}  \otimes \dots \otimes e_{a_m})
 = \sum_{b_1<n_1' ,\dots,b_m <n_m'}
A_{a_1,\dots,a_m}^{b_1,\dots,b_m}\, e_{b_1} \otimes \dots \otimes
e_{j_m}.$$
We will call $(a_1, \dots, a_m, b_1, \dots, b_m)$
{\em the coordinates} of the matrix entry $A_{a_1,\dots,a_m}^{b_1,\dots,b_m}$
of $A$, with respect to the given basis.

The {\em building block} of our construction is a pair of functions
$f_\pm\co \BZ^5 \to \BZ[v^{\pm 1/2}]$, given by
\begin{align*}
f_+(n_1,n_2;&a,b,k) \\&:= (-1)^k v^{-((n_1-1-2a)(n_2-1-2b)+k(k-1))/2}
\qbinom{b+k}{k} \{n_1-1+k-a\}_k, \\
f_- (n_1,n_2;&a,b,k) \\&:=  v^{((n_1-1-2a-2k)(n_2-1-2b+2k)+k(k-1))/2}
\qbinom{a+k}{k} \{n_2-1+k-b\}_k.
\end{align*}
The reader should not focus on the actual, cumbersome formulas.
The main point is that:

{\bf Fact 1}

\begin{itemize}
\item
$f_+$ and $f_-$ are $q$--proper hypergeometric and thus
$q$--holonomic in all variables.
\end{itemize}

For each pair $(n_1,n_2)\in \BN^2$ we define two operators
$$
\BB_+(n_1,n_2),\BB_-(n_1,n_2)\co  V(n_1)\otimes V(n_2) \to V(n_2)
\otimes V(n_1)
$$
by
\begin{eqnarray*}
(\BB_+(n_1,n_2))_{a,b}^{c,d}&:=& f_+(n_1,n_2;a,b,c-b) \, \delta_{c-b,a-d},
\\
(\BB_-(n_1,n_2))_{a,b}^{c,d}&:=& f_+(n_1,n_2;a,b,b-c) \, \delta_{c-b,a-d},
\end{eqnarray*}
where $\delta_{x,y}$ is Kronecker's delta function. Although the coordinates
$(a,b,c,d)$ of the entry $\BB_\pm(n_1,n_2))_{a,b}^{c,d}$ of the operators
$\BB_\pm(n_1,n_2)$ are defined
for $0 \leq a,b \leq n_1$ and $0 \leq c,d \leq n_2$, the above formula makes
sense for all non-negative integers $a,b,c,d$. This will be important for us.
The following lemma is obvious.

\begin{lemma}
\label{braid}
The discrete functions $\BB_\pm(n_1,n_2)^{c,d}_{a,b}$
are $q$--holonomic with respect to the variables $(n_1,n_2,a,b,c,d)$.
\end{lemma}

If we identify $V(n)$ with the simple $n$--dimensional
$U_q(\msl_2)$--module, with $e_i, i=0,\dots,n-1$ being the {\em standard
basis},  then $\BB_+(n_1,n_2),\BB_-(n_1,n_2)$ are
respectively the braiding operator and its inverse acting on $V(n_1)\otimes
V(n_2)$. This fact follows from the formula of the $R$--matrix,
say, in \cite[Chapter 3]{Jantzen}. In particular, $\BB_-(n_1,n_2)$
is the inverse of $\BB_+(n_1,n_2)$. If one allows $a,b,c,d$ in
$\BB_\pm(n_1,n_2)_{a,b}^{c,d}$ to run the set $\BN$, then
$\BB_\pm(n_1,n_2)_{a,b}^{c,d}$ define the braid action on the
Verma module corresponding to $V(n_1), V(n_2)$.

Let $B_m$ be the braid group on $m$ strands, with standard
generators $\sigma_1,...,\sigma_{m-1}$:
$$
\s_i=\psdraw{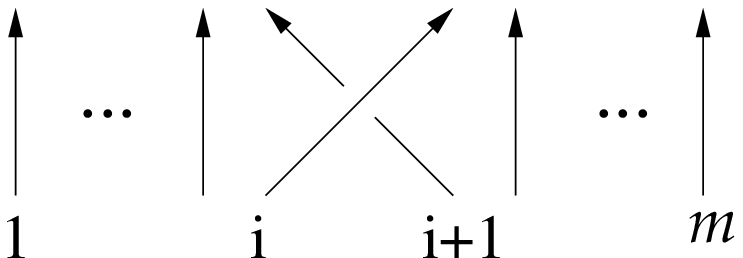}{1.3in}
$$
For each braid $\beta\in
B_m$ and $(n_1,\dots,n_m)\in \BN^m$, we will define an operator
$\tau(\beta)=\tau(\beta)(n_1,\dots,n_m)$,
$$
\tau(\beta) \co  V(n_1)\otimes \dots \otimes V(n_m) \to
V(n_{\bar\beta(1)}) \otimes \dots \otimes V(n_{\bar \beta (m)}),
$$
where $\bar\beta$ is the permutation of $\{1,\dots,m\}$
corresponding to $\beta$. The operator $\tau(\b)$ is uniquely
determined by the following properties: For an elementary braid
$\sigma_i$, we have:
$$
\tau(\sigma_i^{\pm 1}) = \id^{\otimes i-1}\otimes \BB_{\pm }(n_i,n_{i+1})
\id^{\otimes m-i-1}.
$$
In addition, if $\beta=\beta'\beta''$, then $\tau(\beta):=
\tau(\beta')\tau(\beta'')$. It is well-known that $\tau(\beta)$ is
well-defined.

From Fact 0 and Lemma  \ref{braid} it follows that

\begin{lemma}
For any  braid   $\beta\in B_m$,  the discrete function
$\tau(\beta)(n_1,\dots,n_m)$, considered as a function with
variables $n_1,\dots,n_m$ and all the coordinates of the matrix
entry, is $q$--holonomic.
\end{lemma}

Let $K$ be the linear endomorphism of $V(n_1)\otimes \dots
\otimes V(n_m)$ defined by
$$K(e_{i_1}\otimes \dots \otimes
e_{i_m}) = v^{n_1+\dots+n_m -2i_1-\dots -2i_m-m} \, e_{i_1}\otimes
\dots \otimes e_{i_m}.$$
The inverse operator $K^{-1}$ is well-defined.

\begin{corollary}
For any braid $\beta \in B_m$, the  discrete function
$$
\tilde
\tau(\beta):= \tau(\beta)(n_1,\dots,n_m) \times K^{-1}
$$
is $q$--holonomic in $n_1,\dots,n_m$ all all of the coordinates of the
matrix entry.
\end{corollary}

In general, the trace of $\tilde \tau(\beta)$ is called the
{\em quantum trace} of $\tau(\beta)$. Although the target space and
source space maybe different, let us define the quantum trace of
$\tau(\beta)(n_1,\dots,n_m))$ by
$$
\tr_q(\beta)(n_1,\dots,n_m):= \sum_{1\le i \le m}\sum_{0\le a_i<
n_i}\tilde
\tau(\beta)(n_1,\dots,n_m)^{a_1,\dots,a_m}_{a_1,\dots,a_m}.
$$
It follows from Fact 0 that $\tr_q(\beta)(n_1,\dots,n_m)$ is
$q$--holonomic in $n_1,\dots,n_m$. Restricting this function on
the diagonal defined by $n_i= n_{\bar\beta i}, i=1,\dots,m$, we
get a new  function $J_\beta$ of $p$ variables, where $p$ is the
number of cycles of the permutation $\bar \beta$.

Suppose a framed link $L$ can be obtained by closing the braid
$\beta$. Then the colored Jones polynomial $J_L$ is exactly
$J_{\b}$.  Hence Theorem \ref{thm.1} follows.
\end{proof}

\begin{remark}
In general, $J_\KK(n)$ contains the
fractional power $q^{1/4}$. If $K$ has framing 0, then
$J_{\KK'}(n):= J_\KK(n)/[n] \in \BZ[q^{\pm 1}]$. See \cite{Le1}.
\end{remark}

\begin{remark}
\label{rem.break}
There is a variant of the colored Jones function $J_{L'}$ of a colored link
$L'$ where one of the components is broken. If $\b$ is a braid as above,
let us define the {\em broken quantum trace} $\tr'_{\b}$ by
$$
\tr'_q(\beta)(n_1,\dots,n_m):= \sum_{2\le i \le m}\sum_{0\le a_i<
n_i}\tilde
\tau(\beta)(n_1,\dots,n_m)^{a_1,\dots,a_m}_{a_1,\dots,a_m}|_{a_1=0}.
$$
Restricting this function on
the diagonal defined by $n_i= n_{\bar\beta i}, i=1,\dots,m$, we
get a new  function $J_{\beta'}$ of $p$ variables, where $p$ is the
number of cycles of the permutation $\bar \beta$.

If $L'$ denotes the broken link which is the closure of all but the first
strand of $\b$, then the colored Jones function $J_{L'}$ of $L'$
satisfies $J_{L'}=J_{\b'}$.

If $L$ denotes the closure of the broken link $L'$, then we have:
$$
J_L =J_{L'} \times [\l] 
$$
where $\l$ is the color of the broken component of $L'$.
\end{remark}

\subsection{A multisum formula for the colored Jones function of a knot}
\label{sub.multisum}

In this section we will give explicit multisum formulas for the
$\mathfrak{sl}_2$--colored Jones function of a knot. The
calculation here is computerized in Section \ref{sec.computer}.

Consider a word
$w=\s_{i_1}^{\ep_1} \dots \s_{i_c}^{\ep_c}$ of length $m$ written in the standard
generators $\s_1,\dots,\s_{s-1}$ of the braid group $B_m$ with $m$ strands,
where $\ep_i = \pm 1$ for all $i$.

$w$ gives rise to a braid $\b \in B_m$, and
we assume that the closure of $\b$ is a knot $\KK$. Let $\KK'$ denote
{\em long knot} which is the closure of all but the first strand of $\b$.

A {\em coloring} of $\KK'$ is a tuple $\bk=(k_1,\dots,k_c)$ of
{\em angle variables} placed at the crossings of $\KK'$.

\begin{lemma}
\label{lem.1}
There is a unique way to extend a coloring $\bk$ of
$\KK'$ to a coloring of the crossings and part-arcs of $\KK'$
such that
\begin{itemize}
\item
around each crossing the
following consistency relations are satisfied:
$$
\psdraw{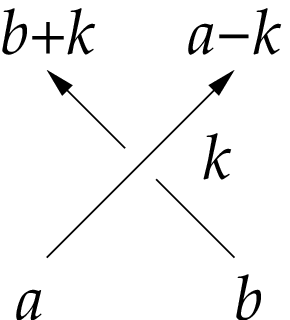}{0.6in} \hspace{1.2cm} \psdraw{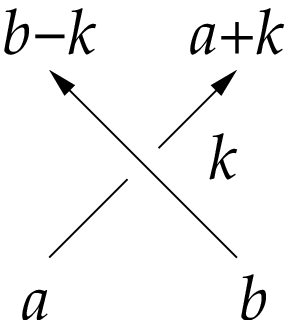}{0.6in}
$$
\item
The color of the lower-left incoming part-arc is $0$.
\end{itemize}
Moreover, the labels of the part-arcs
are linear forms on $\bk$.
\end{lemma}

\begin{proof}
Start walking along the long knot starting at the incoming part-arc. At the
first crossing, whether over or under, the label of the outgoing part-arc is
determined by the label of the ending part-arc and the angle variable of
the crossing.
Thus, we know the label of the outgoing part-arc of the first crossing. Keep
going. Since $\KK'$ is topologically an interval, the result follows.
\end{proof}

\begin{figure}[ht!]\small
$$
w=\s_1^3 \hspace{1cm}
\b=\psdraw{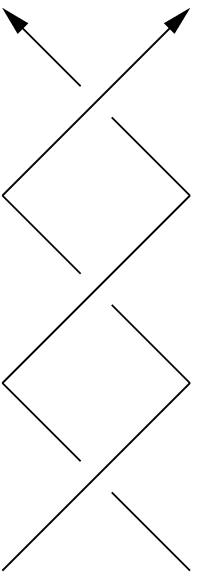}{0.3in} \hspace{1cm}
\KK'=\psdraw{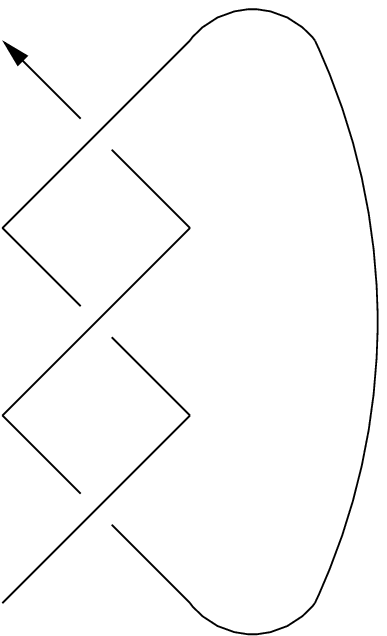}{0.45in} \hspace{1cm}
\psdraw{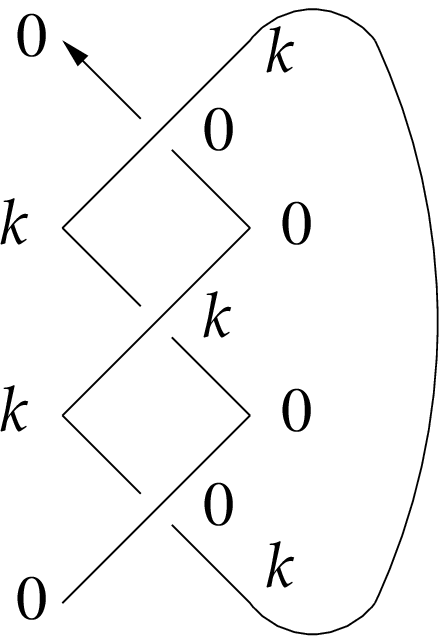}{1in}
$$
\caption{A word $w$, the corresponding braid $\b$, its long closure
$\KK'$, and a coloring of $\KK'$}\label{f.trefoil}
\end{figure}

For an example, see Figure \ref{f.trefoil}.

Fix a coloring of $\KK'$ determined by a vector $\bk$. Let
$b_i(\bk)$ for $i=1,\dots,m$ denote the labels of the top part-arcs of $\b$.
Let $x_j(\bk)$ and $y_j(\bk)$ denote the labeling
of the left and right incoming part-arcs at the $i$th crossing of $\KK'$
for $j=1,\dots,c$. According to Lemma \ref{lem.1}, $b_i(\bk), x_j(\bk)$
and $y_j(\bk)$ are linear forms on $\bk$.

It is easy to see that
$$
\tr'_q(\b)=\sum_{\bk \geq \boldsymbol{0}} F_w(n,\bk)
$$
where
$$
F_w(n,\bk):=
\prod_{i=2}^m v^{\frac{n}{2}-b_i(\bk)}
\prod_{j=1}^c f_{\text{sgn}(\ep_i)}(n,n;x_j(\bk),y_j(\bk)).
$$
is a $q$--proper hypergeometric function.
Remark \ref{rem.break} then implies that

\begin{proposition}
\label{prop.multi}
The colored Jones function of a long knot $\KK'$ is a multisum
of a $q$--proper hypergeometric function:
$$
J_{\KK'}(n)=\sum_{\bk \geq \mathbf{0}} F_w(n,\bk).
$$
\end{proposition}

\begin{remark}
\label{rem.multi}
If a long knot $\KK'$ is presented by a planar projection $D$
with $c$ crossings (which is not necessarily the closure of a braid),
then similar to the above there is a $q$--proper hypergeometric function
$F_D(n,\bk)$ of $c+1$ variables such that
$J_{\KK'}(n)=\sum_{\bk \geq \mathbf{0}} F_D(n,\bk)$.
Of course, $F_D$ depends on the planar projection.
Occasionally, some of the summation variables can be ignored. This is the case
for the right hand-trefoil (where the multisum reduces to a single sum)
and the figure eight (where it reduces to a double sum).
\end{remark}

D Bar-Natan has kindly provided us with a computerized version of Proposition
\ref{prop.multi}, \cite{BN}.

\section{The cyclotomic function of a knot is $q$--holonomic}
\label{sec.cyclo}

Habiro \cite{Habiro} proved that the colored Jones polynomial (of
$\msl_2$) can be rearranged in the following convenient form,
known as the {\em cyclotomic expansion} of the colored Jones polynomial:
For every $0$--framed knot $\KK$, there exists a function
$$
C_\KK\co \BZ_{>0}\to \BZ[q^{\pm1}]
$$
$$J_\KK(n) =\sum_{k=1}^\infty C_\KK(k) S(n,k),
\leqno{\hbox{such that}}$$
$$S(n,k):= \{n+k-1\}_{2k-1}/(v-v^{-1}) =
\frac{\prod_{n-k+1}^{n+k-1} \,(v^i-v^{-i})}{v-v^{-1}}.
\leqno{\hbox{where}}$$
Note that $S(n,k)$ does not depend on the knot $\KK$. Note that
$J$ is determined from $C$ and vice-versa by an upper diagonal
matrix, thus $C$ takes values in $\BQ(q)$. The difficult part of
Habiro's result is  $C_\KK$ takes values in $\BZ[q^{\pm}]$. The
integrality of the cyclotomic function
is a crucial ingredient in the study of
integrality properties of 3--manifold invariants, \cite{Habiro}.

\begin{theorem}
\label{thm.cyclo}
The cyclotomic function $C_\KK\co \BN\to \BZ[q^{\pm}]$ of every knot $\KK$ is
$q$--holonomic.
\end{theorem}

\begin{proof}
Habiro showed that $C_\KK(n)$ is the quantum invariant
of the knot $\KK$ with color
$$
P''(n):= \frac{\prod_{i=1}^{n-1} (V(2)- v^{2i-1}
-v^{1-2i})}{\{2n-1\}_{2n-2}},
$$
where $V(n)$ is the unique $n$--dimensional  simple
$\msl_2$--module, and (retaining Habiro's notation with a shift
$n\to n-1$) $P''(n)$ is considered as an element of the ring of
$\msl_2$--modules over $\BQ(v)$.

Using induction one can easily prove that
$$ P''(n)= \sum_{k=1}^n R(n,k) V(k),$$
where  $R(k,n)$ is given by
$$
R(n,k)= (-1)^{n-k} \frac{\{2k\} }{\{2n-1\}![2n] }\qbinom{2n}{n-k}.
$$
We learned this formula from Habiro \cite{Habiro} and Masbaum \cite{Ma}.
Since
$$
C_\KK(n) =  \sum_{k} R (n,k) J_\KK(k)
$$
and $R(n,k)$ is $q$--proper hypergeometric  and thus
$q$--holonomic in both variables $n$ and $k$, it follows that $C_\KK$ is
$q$--holonomic.
\end{proof}

\section{Complexity}
\label{sec.effective}

In this section we show that Theorem \ref{thm.1} is effective. In other words,
we give a priori bounds and computations that appear in Theorem \ref{thm.c}.

\subsection{Finding a recursion relation for multisums}
\label{sub.u1}

Our starting point are multisums of $q$--proper hypergeometric functions.
Recall the definition \ref{def.proper} of a $q$--proper hypergeometric function
$F(n,\bk)$ from Section \ref{sub.hyper}, and let $G$ denote
$$
G(n):=\sum_{\bk\geq \mathbf{0}} F(n,\bk)
$$
throughout this section.

With the notation of Equation \eqref{eq.proper},
Wilf--Zeilberger show that:

\begin{theorem}
[{\cite[Sec.5.2]{WZ}}]\label{thm.WZ}$\phantom{99}$

{\rm (a)}\qua $F(n,\bk)$ satisfies a $k$--free
recurrence relation of order at most
$$
J^\star:=\frac{(4 S T B^2)^r}{r!}
$$
where $B=\max_{s,t}\{ |\mathbf{b}_s|, |\mathbf{v}_t|, |a_s|, |u_t|\}+
\max_{\mu,\nu} |a_{\mu,\nu}|$ where $a_{\mu,\nu}$ are the coefficients of the
quadratic form $A$.

{\rm (b)}\qua
Moreover, $G(n)$ satisfies an inhomogeneous
recursion relation of order at most $J^\star$.
\end{theorem}

Let us briefly comment on the proof of this theorem.
A {\em certificate} is an operator of the form
$$
P(E,Q) + \sum_{i=1}^r (E_i-1) R_i(E,E_1,\dots,E_r,Q,Q_1,\dots,Q_r)
$$
that annihilates $F(n,\bk)$, where $P$ and $R_i$ are operators
with $P$ a polynomial in $E,Q$, with
$P \neq 0$. Here $E$ is the shift operator on $n$, $E_i$ (for $i=1,
\dots,r$) are shift operators in $k_i$, and $Q$ is the multiplication
operator by $q^{k}$ and  $Q_i$ (for $i=1, \dots,r$) are
the multiplication operator by $q^{k_i}$, where $\bk=(k_1,\dots,k_r)$.

The important thing is that $P(E,Q)$ is an operator that does not depend on
the summation variables $\bk$. A certificate implies that for all $(n,\bk)$
we have:
\begin{align*}
P(E,Q) F (n,\bk) + \sum_{i=1}^r (
G_i(n,k_1,\dots,k_{i-1}, k_i+1,k_{i+1},\dots, k_r)\,\,-&\\
G_i(n,k_1,\dots,k_{i-1}, k_i,k_{i+1},\dots, k_r) )&=0,
\end{align*}
where $G_i(n,\bk)=R_i F(n,\bk)$. Summing over $\bk \geq
\boldsymbol{0}$, it follows that $ G(n) $ satisfies an
inhomogeneous recursion relation $PG=\text{error}(n)$. Here
$\text{error}(n)$ is a sum of multisums of $q$--proper
hypergeometric functions of one variables less. Iterating the
process, we finally arrive at a homogeneous recursion relation for
$G$.

How can one find a certificate given $F(n,\bk)$? Suppose that $F$
satisfies a $\bk$--{\em free} recursion relation $AF=0$, where
$A=A(E,Q,E_1,\dots,E_r)$ is an operator that does not depend on
the $Q_i$. Then, evaluating $A$ at $E_1=\dots E_r=1$, we obtain
that
$$
A=A(E,Q,1,\dots,1) + \sum_{i=1}^r (E_i-1) R_i(E,Q,E_1,\dots,E_r)
$$
is a certificate.

How can we find a $\bk$--free recursion relation for $F$? Let us write
$$
A=\sum_{(i,\bs{j}) \in S} \s_{i,\bs{j}}(Q) E^i E^{\bs{j}}
$$
where $S$ is a finite set, $\bs{j}=(j_1,\dots,j_r)$, $E^{\bs{j}}=
E_1^{j_1} \dots E_r^{j_r}$, and $\s_{i,\bs{j}}(Q)$ are polynomial functions
in $Q$ with coefficients in $\BQ(q)$; see \cite{R}.
The condition $AF=0$ is equivalent to the equation
$(AF)/F=0$. Since $F$ is $q$--proper
hypergeometric, the latter equation is the vanishing of a rational function
in $Q_1,\dots,Q_r$. By cleaning out denominators, this is
equivalent to a system of {\em linear
equations} (namely, the coefficients of monomials in $Q_i$ are zero),
with unknowns the polynomial functions $\s_{i,\bs{j}}$. For a careful
discussion, see \cite{R}.
As long as there are more unknowns than equations, the system
is guaranteed to have a solution. \cite{WZ} estimate the number of equations
and unknowns in terms of $F(n,\bk)$, and prove Theorem \ref{thm.WZ}.

Wilf--Zeilberger programmed the above proof, see \cite{PWZ}. As
time passes the algorithms get faster and more
refined. For the state-of-the-art algorithms and implementations, see
\cite{PR1,PR2} and \cite{R}, which we will use below.

Alternative algorithms of noncommutative elimination, using {\em noncommutative
Gr\"obner basis}, have been developed by Chyzak and Salvy, \cite{CS}.
In order for have Gr\"obner basis, one needs to use the following localization
of the $q$--Weyl algebra
$$
\B_r= \frac{Q(q,Q_1, \dots, Q_r)\la E_1, \dots, E_r
\ra}{ (\text{Rel}_q)}.
$$
and {\em Gr\"obner basis} \cite{CS}.

In case $r=1$, $\B_1$ is a {\em principal ideal domain} \cite[Chapter 2,
Exercise 4.5]{Cou}. In that case one can associate an operator in $\B_1$
(unique up to units) that generates that annihilating ideal of $G(n)$.
For a conjectural relation between this operator for the $\msl_2$--colored
Jones function of a knot and hyperbolic geometry, see \cite{G2}.

Let us point out however that none of the above algorithms can find generators
for the annihilating ideal of the multisum $G(n)$. In fact, it is an open
problem
how to find generators for the annihilating ideal of $G(n)$ in terms of
generators for the annihilating ideal of $F(n,\bk)$,
in theory or in practice. We thank M Kashiwara for pointing this out to us.

\subsection{Upper bounds for initial conditions}
\label{sub.u2}

In another direction, one may ask the following question: if a $q$--holonomic
function satisfies a nontrivial recursion relation, it follows that it
is uniquely determined by a finite number of initial conditions. How many?
This was answered by Yen, \cite{Y}. If $G$ is a discrete
function which satisfies a recursion relation of order $J^\star$, consider
its principal symbol $\s(q,Q)$, that is the coefficient of the leading
$E$--term.
The principal symbol lies in the commutative ring $\BZ[q^{\pm},Q^{\pm}]$
of Laurrent polynomials in two variables $q$ and $Q$. For every $n$, consider
the Laurrent polynomial $\s(q,q^n) \in \BZ[q^{\pm}]$. If $\s(q,q^n) \neq 0$
for all $n$, then $G$ is determined by $J^\star$ many initial
values. Since $\s(q,Q) \neq 0$, it follows that $\s(q,q^n) \neq 0$
for large enough $n$. In fact, in \cite[Prop.3.1]{Y} Yen proves that
$\s(q,q^n) \neq 0$ if
$n > \deg_q(\s)$, then $\s(q,q^n) \neq 0$, where the degree of a Laurrent
polynomial in $q$ is the difference between the largest and smallest exponent.
Thus, $G$ is determined by $J^{\star\star}:=J^\star + \deg_q(\s)$
initial conditions.

Yen further gives upper bounds for $\deg_q(\s)$ in terms of the
$q$--hypergeometric summand, see \cite[Thm.2.9]{Y} for single sums. 
An extension of Yen's work 
to multisums, gives a priori upper bounds $J^{\star\star}$
in terms of the $q$--hypergeometric summand. These exponential
bounds are of theoretical interest only, and in practice much smaller
bounds are found by computer.

\subsection{Proof of Theorem \ref{thm.c}}
\label{sub.thm.c}

Theorem \ref{thm.c} follows from Proposition \ref{prop.multi} together
with the discussion of Sections \ref{sub.u1} and \ref{sub.u2}.
\endproof

Our luck with the colored Jones function is that we can identify
it with a multisum of a $q$--proper hypergeometric function. Are we
really lucky, or is there some deeper explanation? We believe that
there is a underlying geometric reason for coincidence, which in a
sense explains the underlying geometry of topological quantum
field theory. We will postpone to a later publication applications
of this principle to Hyperbolic Geometry; \cite{G2}.

\section{In computer talk}
\label{sec.computer}

In this section we will show that Proposition \ref{prop.multi} can be
implemented by computer.

For every knot, one can write down a multisum formula for the colored Jones
function, where the summand is $q$--hypergeometric.
Occasionally, this multisum formula can be written as a single
sum. There are various programs that can compute the recursion relations
and their orders for multisums. In maple, one may use {\tt qEKHAD} developed by
Zeilberger \cite{PWZ}. In Mathematica, one may use the {\tt qZeil.m} and
{\tt qMultiSum.m} packages of RISC developed by Paule and Riese
\cite{PR1,PR2,R}.

\subsection{Recursion relations for the cyclotomic function of twist knots}
\label{sub.cyclo}

\begin{figure}[ht!]
$$
\psdraw{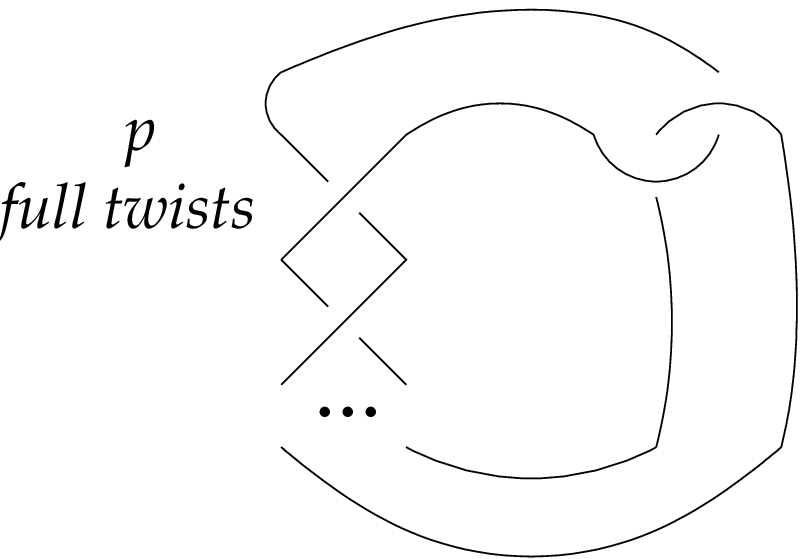}{2in}
$$
\caption{The twist knot $K_p$, for integers $p$. For $p=-1$, it is the
Figure 8, for $p=0$ it is the unknot, for $p=1$ it is the left trefoil
and for $p=2$ it is the Stevedore's ribbon knot.}\label{twist}
\end{figure}

The twist knots $Kp$ for integer $p$ are shown in Figure
\ref{twist}. Their planar projections have $2|p|+2$ crossings,
$2|p|$ of which come from the full twists, and $2$ come from the
negative clasp.

Masbaum, \cite{Ma}, following Habiro and Le gives the
following formula for the cyclotomic function of a twist knot. Let
$c(p, \cdot)$ denote the cyclotomic function of the twist knot $Kp$.
Rearranging a bit Masbaum's formula \cite[Eqn.(35)]{Ma}, we obtain that:
\begin{align}
c(p,n) &=(-1)^{n+1} q^{n(n+3)/2} \notag\\
&\sum_{k=0}^\infty (-1)^k q^{k(k+1)p+k(k-1)/2}
(q^{2k+1}-1)\frac{(q;q)_n}{(q;q)_{n+k+1} (q;q)_{n-k}}
\label{eq.twist}
\end{align}
The above sum has compact support for each $n$. Now, in computer talk,
we have:

{\small
\begin{verbatim}
Mathematica 4.2 for Sun Solaris
Copyright 1988-2000 Wolfram Research, Inc.
 -- Motif graphics initialized --
In[1]:=<< qZeil.m
\end{verbatim}
\verb+q-Zeilberger Package by Axel Riese --+ \copyright \verb+RISC Linz -- V 2.35 (04/29/03)+
}

\vskip 2mm
For $p=-1$ (which corresponds to the Figure 8 knot) the program gives:
\vskip 2mm

{\small
\begin{verbatim}
In[2]:= qZeil[q^(n(n + 3)/2) (-1)^(n + k + 1) q^(-k(k + 1))(q^(2k + 1) 
        - 1)qfac[q, q, n]/(qfac[q, q, n + k + 1] qfac[q, q, n - k]) 
        q^(k(k - 1)/2), {k, 0, Infinity}, n, 1]

Out[2]= SUM[n] == SUM[-1 + n]
\end{verbatim}
}

which means that $c(-1,n)=c(-1,n-1)$ in accordance to the discussion
after \cite[Thm.5.1]{Ma} which states $c(-1,n)=1$ for all $n$.

\vskip 2mm
For $p=1$ (which corresponds to the left hand trefoil) the program gives:
\vskip 2mm

{\small
\begin{verbatim}
In[3]:= qZeil[q^(n(n + 3)/2) (-1)^(n + k + 1) q^(k(k + 1))(q^(2k + 1) 
        - 1)qfac[q, q, n]/(qfac[q, q, n + k + 1] qfac[q, q, n - k]) 
        q^(k(k - 1)/2), {k, 0, Infinity}, n, 1]

                     1 + n
Out[3]= SUM[n] == -(q      SUM[-1 + n])
\end{verbatim}
}

which means that $c(1,n)=-q^{n+1}c(1,n-1)$ in accordance to the
discussion after \cite[Thm.5.1]{Ma} which states $c(1,n)=(-1)^n
q^{n(n+3)/2}$ for all $n$.

Similarly, for $p=2$ (which corresponds to Stevedore's ribbon knot) the program
gives:

{\small
\begin{verbatim}
In[4]:= qZeil[q^(n(n + 3)/2) (-1)^(n + k + 1) q^(2k(k + 1))(q^(2k + 1) 
        - 1)qfac[q, q, n]/(qfac[q, q, n + k + 1] qfac[q, q, n - k]) 
        q^(k(k - 1)/2), {k, 0, Infinity}, n, 1]
Out[4]:= No solution: Increase order by 1
\end{verbatim}
}

which {\em proves} that $c(2,n)$ satisfies no first order recursion relation.
It does satisfy a second order recursion relation, as we find by:

{\small
\begin{verbatim}
In[5]:= qZeil[q^(n(n + 3)/2) (-1)^(n + k + 1) q^(2k(k + 1))(q^(2k + 1) 
        - 1) qfac[q, q, n]/(qfac[q, q, n + k + 1] qfac[q, q, n - k]) 
        q^(k(k - 1)/2), {k, 0, Infinity}, n, 2]
                     2 + 2 n       -1 + n
Out[4]= SUM[n] == -(q        (1 - q      ) SUM[-2 + n]) -

       1 + n           n    2 n
>     q      (1 + q - q  + q   ) SUM[-1 + n]
\end{verbatim}
}

Thus, the program computes not only a recursion relation, but also the
order of a minimal
one. Experimentally, it follows that $c(p,n)$ satisfies a recursion
relation of order $|p|$, for all $p$. Perhaps one can guess the form
of a minimal order recursion relation for all twist knots.

Actually, more is true. Namely, the formula for $c(p,n)$ shows that
it is a $q$--holonomic function in {\em both} variables $(p,n)$. Thus, we are
guaranteed to find recursion relations with respect to $n$ and with
respect to $p$. Usually, recursion relations with respect to $p$ for fixed
$n$ are called {\em skein theory} for the $n$th colored Jones function,
because the knot is changing, and the color is fixed.

Thus, $q$--holonomicity implies skein relations (with respect to the number
of twists) for the $n$th colored Jones polynomial of twist knots, for
every fixed $n$.

For computations of recursion relations of the cyclotomic function
of twist knots, we refer the reader to \cite{GaS}.

\subsection{Recursion relations for the colored Jones function of the figure
8 knot}
\label{sub.figure8}

The Mathematica package {\tt qMultiSum.m} can compute recursion relations for
$q$--multisums.
Using this, we can compute equally easily the recursion relation for the
colored Jones function. Due to the length of the output, we illustrate this by
computing the recursion relation for the colored Jones function of the Figure
8 knot. Recall from Equation \eqref{eq.twist} for $p=-1$ and from the fact
that $c(-1,n)=1$ that the colored Jones function of the figure 8 knot is given
by:
\begin{equation}
\label{eq.figure8}
J_{K(-1)}(n)=\sum_{k=0}^\infty q^{nk} (q^{-n-1};q^{-1})_k (q^{-n+1};q)_k.
\end{equation}
In computer talk,

{\small
\begin{verbatim}
In[6]:= qZeil[q^(n k) qfac[q^(-n-1),q^(-1),k] qfac[q^(-n+1),q,k],
        {k,0,Infinity},n,2]

                   -1 - n       n         2 n
                  q       (q + q ) (-q + q   )
Out[6]= SUM[n] == ---------------------------- -
                                  n
                            -1 + q

            -2 + n        -1 + 2 n
      (1 - q      ) (1 - q        ) SUM[-2 + n]
>     ----------------------------------------- +
                    n        -3 + 2 n
              (1 - q ) (1 - q        )

        -2 - 2 n       -1 + n 2       -1 + n
>     (q         (1 - q      )  (1 + q      )
\end{verbatim}
\begin{verbatim}
      4    4 n    3 + n    1 + 2 n    3 + 2 n    1 + 3 n
>   (q  + q    - q      - q        - q        - q       ) SUM[-1 + n])/

              n        -3 + 2 n
>      ((1 - q ) (1 - q        ))
\end{verbatim}
}

This is a second order {\em inhomogeneous} recursion relation for the colored
Jones function. A third order homogeneous relation may be obtained by:

{\small
\begin{verbatim}
In[7]:= MakeHomRec[%, SUM[n]]

         2 + n    3    n
        q      (-q  + q ) SUM[-3 + n]
Out[7]= ----------------------------- -
             2    n     5    2 n
           (q  + q ) (-q  + q   )

        -2 - n   2    n    8    4 n      6 + n    7 + n    3 + 2 n
>     (q       (q  - q ) (q  + q    - 2 q      + q      - q        +

            4 + 2 n    5 + 2 n    1 + 3 n      2 + 3 n
>          q        - q        + q        - 2 q       ) SUM[-2 + n]) /

              n    5    2 n
>      ((q + q ) (q  - q   )) +

        -1 - n        n    4    4 n    2 + n      3 + n    1 + 2 n
>     (q       (-q + q ) (q  + q    + q      - 2 q      - q        +

            2 + 2 n    3 + 2 n      1 + 3 n    2 + 3 n
>          q        - q        - 2 q        + q       ) SUM[-1 + n]) /

                                  1 + n        n
          2    n         2 n     q      (-1 + q ) SUM[n]
>      ((q  + q ) (-q + q   )) + ----------------------- == 0
                                         n        2 n
                                   (q + q ) (q - q   )
\end{verbatim}
}
Of course, we can clear denominators and write the above recursion relation
using the $q$--Weyl algebra $\A$. Let us end with a matching the theoretical
bound for the recursion relation from Section \ref{sec.effective} with the
computer calculated bound from this section. Using Theorem \ref{thm.WZ}, it
follows that the summand satisfies a recursion relation of order
$J^\star=1^2+1^2=2$. This implies that the colored Jones function of the
Figure 8 knot satisfies an inhomogeneous relation of degree $2$ as was found
above. The program also confirms that the colored Jones function of the Figure
8 knot does not satisfy an inhomogeneous relation of order less than $2$.

\section{The colored Jones function for a simple Lie algebra}
\label{sec.lie}

Fix a {\em simple complex Lie algebra} $\fg$ of {\em rank} $\ell$. For every
knot $\KK$ and every finite-dimensional $\fg$--module $V$, called
{\em the color} of the knot, one can define the {\em quantum invariant}
$J_\KK(V)\in \BZ[q^{\pm 1/2D}]$, where $D$ is the determinant of
the Cartan matrix of $\fg$. Simple $\fg$--modules are parametrized
by the set of {\em dominant weights}, which can be identified, after we choose
fixed fundamental weights, with $\BN^\ell$. Hence $J_\KK$ can be
considered as a function $J_\KK\co \BN^\ell \to \BZ[q^{\pm 1/2D}]$.

\begin{theorem}
\label{mainn} For every simple Lie algebra other than $G_2$, and a
set of fixed fundamental weights, the colored Jones function
$J_\KK\co  \BN^\ell \to \BZ[q^{\pm 1/2D}]$ is $q$--holonomic.
\end{theorem}

Hence the colored Jones function will satisfy some recursion
relations, which, together with values at a finitely many initial
colors, totally determine the colored Jones function $J_\KK$.

\begin{remark}
\label{rem.G2}
The reason we exclude the $G_2$ Lie algebra is technical. Namely,
at present we cannot prove that the structure constants of
the multiplication of the quantized enveloping algebra of $G_2$
with respect to a standard PBW basis, are $q$--holonomic; see
Remark \ref{rem.appG2}.
We believe however, that the theorem also holds for $G_2$.
\end{remark}

The  proof occupies the rest of this section.  We will define
$J_\KK$ using representation of the braid groups coming from the
$R$--matrix acting on Verma modules (instead of finite-dimensional
modules). We then show that the $R$--matrix is $q$--holonomic. The
theorem follows from that fact that  products and traces of
$q$--holonomic matrices are $q$--holonomic.

\subsection{Preliminaries}
\label{sub.preliminaries}

Fix a {\em Cartan subalgebra} $\hh$ of $\fg$  and a basis
$\{\al_1,\dots,\al_\ell\}$ of simple roots for the dual space
$\hh^*$. Let $\hh^*_\BR$ be the $\BR$--vector space spanned by
$\al_1,\dots,\al_\ell$. The {\em root lattice} $Y$ is the
$\Z$--lattice generated by $\{\al_1,\dots,\a_{\ell}\}$. Let $X$ be
the {\em weight lattice} that is spanned by the {\em fundamental
weights} $\lambda_1,\dots,\lambda_\ell$. Normalize the {\em
invariant scalar product} $(\cdot,\cdot)$ on $\hh^*_\BR$ so that
$(\alpha,\alpha)=2$ for every {\em short} root $\alpha$. Let $D$
be the determinant of the Cartan matrix,  then $(x,y)\in
\frac{1}{D}\BZ$ for $x, y\in X$.

Let $s_i, i=1,\dots,\ell$, be the {\em reflection along the wall}
$\alpha_i^\perp$. The {\em Weyl group} $\Weyl$ is generated by
$s_i,i=1\dots,\ell$, with the braid relations together with $s_i^2=1$. A
word $w= s_{i_1}\dots s_{i_r}$ is {\em reduced} if $w$, considered
as an element of $W$, can not be expressed by a shorter word. In
this case the {\em length} $l(w)$ of the element $w\in W$ is $r$.
The {\em longest element} $\omega_0$ in $W$ has length
$t=(\dim(\fg)-\ell)/2$, the number of positive roots of $\fg$.

\subsubsection{The quantum group  $\UU$}
\label{subsub.UU}

The quantum group
$\UU=\UU_q(\fg)$ associated to $\fg$ is a {\em Hopf algebra} defined
over $\Q(v)$, where $v$ is the usual quantum parameter (see
\cite{Jantzen,Lusztig}). Here our $v$ is the same as $v$ of
Lusztig \cite{Lusztig} and is equal to $q$ of Jantzen
\cite{Jantzen}, while our $q$ is $v^2$. The standard generators of
$\UU$ are $E_\alpha, F_\alpha, K_\alpha$ for $ \alpha \in
\{\alpha_1,\dots,\alpha_\ell\}$. For a full set of relations, as
well as a good introduction to quantum groups, see \cite{Jantzen}.
Note that all the $K_\alpha$'s commute with each other.

For an element $\gamma\in Y$, $\gamma=k_1\alpha_1+\dots+
k_\ell\alpha_\ell$, let $K_\gamma:= K_{\alpha_1}^{k_1}\dots
K_{\alpha_1}^{k_1}$.

There is a $Y$--{\em grading} on $\UU$ defined by $|E_\alpha|=\alpha,
|F_\alpha|=-\alpha$, and $|K_\alpha|=0$. If $x$ is homogeneous,
then
$$K_\gamma x = v^{(\alpha, |x|)} x K_\gamma.$$
Let $\UU^+$ be the subalgebra of $\UU$ generated by the
$E_\alpha$, $\UU^-$ by the $F_\alpha$, and $\UU^0$ by the
$K_\alpha$. It is known that the map
$$\UU^-\otimes \UU^0 \otimes \UU^+ \quad \to \quad \UU$$
$$(x, x', x'') \to x x' x''$$
is an isomorphism of {\em vector spaces}.

\subsubsection{Verma modules and finite dimensional modules}
\label{subsub.verma}

Let $\lambda\in X$ be a weight. The {\em Verma module} $M(\lambda)$
is a $\UU$--module with underlying vector space $\UU^-$ and
with the action of $\UU$ that is uniquely
determined by the following condition.  Here $\eta$ is the unit of the
algebra $\UU^-$:
\begin{eqnarray*}
E_\alpha\cdot \eta &=& 0  \qquad \text{for all} \quad \alpha \\
K_\alpha \cdot \eta &=& v^{(\alpha,\lambda)} \eta \qquad
 \text{for all} \quad \alpha \\
F_\alpha \cdot x &=& F_\alpha x \qquad
\text{for all} \quad \alpha\in \{\alpha_1,\dots,\alpha_\ell\}, x\in \UU^-
\end{eqnarray*}
If $(\lambda|\alpha_i) <0$ for all $i=1,...,\ell$ then
$M(\lambda)$ is irreducible. On the other hand if
$(\lambda|\alpha_i) \ge 0$ for all $i=1,...,\ell$ (ie, $\lambda$
is dominant), then $M(\lambda)$ has a unique proper maximal
submodule, and the quotient $L(\lambda)$ of $ M(\lambda)$ by the
proper maximal submodule is a finite dimensional module (of type
1, see \cite{Jantzen}). Every finite dimensional module of type 1
of $\UU$ is a direct sum of several $L(\lambda)$.

\subsubsection{Quantum braid group action}
\label{subsub.braid}

For each fundamental
root $\alpha \in \{\alpha_1,\dots,\alpha_\ell\}$ there is an {\em
algebra automorphism} $T_\alpha \co  \UU \to \UU$, as described in
\cite[Chapter 8]{Jantzen}. These automorphisms satisfy the
following relations, known as the {\em braid relations}, or {\em Coxeter
moves}.

If $(\alpha,\beta)=0$, then
$T_\alpha T_\beta =T_\beta T_\alpha$.

If $(\alpha,\beta)=-1$, then $T_\alpha T_\beta T_\alpha =T_\beta
T_\alpha T_\beta$.

If $(\alpha,\beta)=-2$, then $T_\alpha T_\beta T_\alpha
 T_\beta=T_\beta T_\alpha T_\beta T_\alpha$.

If $(\alpha,\beta)=-3$, then $
T_\alpha T_\beta T_\alpha T_\beta T_\alpha T_\beta
= T_\beta T_\alpha T_\beta T_\alpha T_\beta T_\alpha $.

Note that the Weyl group is generated by $s_\alpha$ with exactly
the above relations, replacing $T_\alpha$ by $s_\alpha$, and the
extra relations $s_\alpha^2=1$.

Suppose  $w=s_{i_1}\dots s_{i_r}$ is a reduced word, one can
define
$$T_w:= T_{\alpha_{i_1}}\dots T_{\alpha_{i_r}}.$$
Then $T_w$ is well-defined: If $w,w'$ are two reduced words of
the same element in $W$, then $T_w=T_{w'}$. This follows from the
fact that any two reduced presentations of an element of $W$ are
related by a sequence of Coxeter moves.

\subsubsection{Ordering of the roots}
\label{subsub.ordering}

Suppose $w= s_{i_1} s_{i_2} \dots s_{i_t}$ is a reduced word
representing  the longest element $\omega_0$ of the Weyl group.
For $r$ between $1$ and $t$ let
$$
\gamma_r(w):= s_{i_1} s_{i_2} \dots s_{i_{r-1}}(\alpha_{i_r}).
$$
Then the set $\{\gamma_i, i=1,...,t\}$ is exactly the set of
positive roots. We {\em totally order} the set of positive roots
by $\ga_1 < \ga_2 < \dots < \ga_t$. This order depends
on the reduced word $w$, and has the following {\em convexity}
property: If $\beta_1, \beta_2$ are two positive roots such that
$\beta_1+\beta_2$ is also a root, then $\beta_1+\beta_2$ is
between $\beta_1$ and $\beta_2$. In particular, the first and the
last, $\gamma_1$ and $\gamma_t$, are always fundamental roots.
Conversely, any convex total ordering of the set of positive roots comes
from a reduced word representing the longest element of $W$.

\subsubsection{PBW basis for $\UU^-, \UU^+$, and $\UU$}
\label{subsub.PBW}

Suppose
$w=s_{i_1} \dots s_{i_t}$ is a reduced word representing the
longest element of $W$.  Let us define
\begin{eqnarray*}
e_{r}(w) &=& T_{\alpha_{i_1}} T_{\alpha_{i_2}} \dots
T_{\alpha_{i_{r-1}}}(E_{\alpha_{i_r}}), \\
f_{r}(w) &=& T_{\alpha_{i_1}} T_{\alpha_{i_2}} \dots
T_{\alpha_{i_{r-1}}}(F_{\alpha_{i_r}}).
\end{eqnarray*}
Then $|e_r|= \gamma_r=-|f_r|$. (We drop $w$ if there is no
confusion.)

If $\gamma_r$ is one of the fundamental roots,
$\gamma_r=\alpha\in \{ \alpha_1,\dots,\alpha_\ell\}$, then
$e_r(w)=E_\alpha$, $f_r(w)= F_\alpha$ (and do not depend on $w$).

For $t\ge j \ge i \ge 1$ let $\UU^-[j,i]$ be the vector space
spanned by $f_j^{n_j}f_{j-1}^{n_{j-1}}\dots f_i^{n_i}$,
for all $n_j, n_{j-1}, \dots, n_i \in \N$
and let $\UU^+[i,j]$ the vector space  spanned by\break
$e_i^{n_i}e_{i+1}^{n_{i+1}}\dots e_{j}^{n_j}$, for all $n_j,
n_{j-1}, \dots, n_i \in \N$.
It is known that $\UU^- = \UU^-[t,1]$ and $\UU^+=\UU^+[1,t]$.

For $\n=(n_1,\dots,n_t)\in \N^t$, $\j=(j_1,\dots j_\ell) \in
\Z^\ell$ and $\m=(m_1, \dots, m_t) \in \N^t$ let us define
$\f^{\n}$, $\K^{\j}$ and $\e^{\m}$ by
$$
\f^\n(w):= f^{n_t}_t\dots f^{n_1}_1, \qquad
\K^\j:=
K_{j_1\alpha_1}\dots K_{j_\ell \alpha_\ell} \qquad
\e^\n(w):= e^{n_1}_1\dots e^{n_t}_t.
$$
Then as vector spaces over $\BQ(v)$$\UU^-$, $\UU^+$ and $\UU$ have
{\em Poincare-Birkhoff-Witt} (in short, PBW) basis
$$
\{\f^{\n} \, | \,  \n\in \N^t\}, \qquad \{\e^{\m} \, | \, \m\in
\N^t\},  \qquad \{\f^\n \K^\j \e^\m \,\, | \,\, \n,\m\in \N^t,
\j\in \Z^\ell\}
$$
respectively, associated with the reduced word $w$.

In order to simplify notation, we define
$S:= \N^t\times \Z^\ell \times \N^t$, and $x_\s:=\f^n\, K^\j\, \e^m$.
Thus,
\begin{equation}
\label{eq.PBWU}
\{x_\s \, | \, \s \in S\}
\end{equation}
is a PBW basis of $\UU$ with respect to the reduced word $w$.

\subsubsection{A commutation rule}
\label{subsub.commutation}

For $x,y\in \UU$ homogeneous let us define
$$
[x,y]_q:= xy - v^{(|x|,|y|)}yx.
$$
Note that, in general, $[y,x]_q$ is not proportional to $[x,y]_q$.

An important property of the PBW basis is the following
commutation rule, see \cite{Soibelman}.
 If $i<j$ then $[f_i,f_j]_q$
belongs to $\UU^-[j-1,i+1]$ (which is 0 if $j=i+1$). It follows
that $\UU^-[j,i]$ is an {\em algebra}. This
allows us to sort algorithmically non-commutative monomials in the variables
$f_i$. Also two consecutive variables always $q$--commute:
$[f_i,f_{i+1}]_q=0$.

Similarly, if $i<j$ then $[e_i,e_j]_q$ belongs to $\UU^+[i+1,j-1]$
(which is 0 if $j=i+1$). It follows that $\UU^+[i,j]$ is an {\em
algebra}, and two consecutive variables always $q$--commute,
$[e_i,e_{i+1}]_q=0$.

\subsection{$q$--holonomicity of quantum groups}
\label{sub.holo}

Suppose $A\co  \UU\to \UU$ is a linear operator. Using the PBW basis
of $\UU$ (see Equation \eqref{eq.PBWU}),
we can present $A$ by a matrix:
$$
A(x_\s) = \sum_{\s'} A_{\s}^{\s'} x_{\s'},
$$
with $A_{\s}^{\s'} \in \Q(v)$. We call $(\s,\s')$ the coordinates
of the matrix entry $A_{\s}^{\s'}$.

\begin{definition}
\label{def.holo}
We say that $A$ is $q$--{\em holonomic} if the matrix entry
$A_{\s}^{\s'}$, considered as a function of $(\s,\s')$ is
$q$--holonomic with respect to all the variables.
\end{definition}

A priori this definition depends on the reduced word $w$. But we will soon see
that if $A$ is $q$--holonomic in one PBW basis, then it is so in any
other PBW basis.

\subsubsection{$q$--holonomicity of transition matrix}
\label{subsub.transition}

Suppose $x_{\s}(w')$ is another PBW basis associate to another
reduced word $w'$ representing the longest element of $W$. Then
we have the transition matrix $M_{\s}^{\s'}$ between the two
bases, with entries in $\Q(v)$. The next proposition checks that the entries
of the transition matrix are $q$--holonomic, by a standard reduction to the
rank $2$ case.

\begin{proposition}
\label{prop.trans}
Except for the case of $G_2$, the  matrix entry $M_{\s}^{\s'}$ is $q$--holonomic
with respect to all its coordinates.
\end{proposition}
\begin{proof}
Since any two reduced presentations of an element of $W$ are related by a
sequence of Coxeter moves, it is enough to consider the case of a single
Coxeter move. Since each Coxeter move involves only two fundamental roots
and all $T_\alpha$'s are {\em algebra}
isomorphisms, it is enough to considered the case of rank 2 Lie
algebras. For all rank 2 Lie algebras (except $G_2$) we present the proof in
Appendix.
\end{proof}

\subsubsection{Structure constants}
\label{subsub.constants}

Recall the PBW basis $\{x_\s \, | \, \s \in S\}$ of the algebra $\UU$.
The multiplication in $\UU$ is determined by the structure constants
$c(\s,\s',\s'') \in \BQ(v)$ defined by:
$$
x_\s x_{\s'} = \sum_{\s''} c(\s,\s',\s'') x_{\s''}.$$
We will show the following:

\begin{theorem}
\label{main}
The structure constant $c(\s,\s',\s'')$ is $q$--holonomic
with respect to all its variables.
\end{theorem}

Proof will be given in  subsection \ref{proof}.

\subsubsection{Actions on Verma modules are $q$--holonomic}
\label{subsub.actionsverma}

Each Verma module $M(\lambda)$ is naturally isomorphic to $\UU^-$,
as a vector space, via the map $u\to u\cdot \eta$. Using this
isomorphism we identify a PBW basis of $\UU^-$ with a basis  of
$M(\lambda)$, also called a PBW basis. If $u\in \UU$, then the
action of $u$ on $M(\lambda)$ in a PBW basis  can be written by a
matrix $u_\n^{\n'}$ with entries in $\Q(v)$. We call $(\n,\n') \in
\N^t\times \N^t$ the coordinates of the matrix entry.

\begin{proposition}
\label{Verma}
For every $r$ with $1 \leq r \leq t$,
the entries of the matrices $e_r^k, f_r^k$ are $q$--holonomic with
respect to $k, \lambda$, and the coordinates of the entry.
\end{proposition}

This Proposition follows immediately from Theorem \ref{main} and
Fact 0.

\subsection{Quantum knot invariants}
\label{sub.qknots}

\subsubsection{The quasi-$R$--matrix}
\label{subsub.Rmatrix}

Fix a reduced word $w$
representing the longest element of $W$. For each $r, 1\le r \le
t$, let
$$\Theta_r:= \sum_{k\in \N} c_k \, f_r^k \otimes
e_r^k,$$
$$c_k= (-1)^k v_{\gamma_r}^{-k(k-1)/2}\,
\frac{(v_{\gamma_r}-v_{\gamma_r}^{-1})^k} {[k]_{\gamma_r}!}.\leqno{\hbox{where}}
$$
Here $v_{\gamma} = v^{(\gamma|\gamma)/2}$, and
$$ [k]_\gamma! =\prod_{i=1}^k
\frac{v_\gamma^i-v_\gamma^{-i}}{v_\gamma-v_\gamma^{-1 }}.$$
The main thing to observe is that $c_k$ is $q$--holonomic with respect
to $k$. Note
that although $\Theta_r$ is an infinite sum, for every weight
$\lambda\in X$, the action of $\Theta_r$ on  $M(\lambda)\otimes
M(\lambda)$ is well-defined. This is because the action of $e_r$
is locally nilpotent, ie, for every $x\in M(\lambda)$, there is
$k$ such that $e_r^k\cdot x=0$.

The {\em quasi-$R$--matrix} is:
$$\Theta := \Theta_{t} \Theta_{{t-1}} \dots
\Theta_{1}.$$
We will consider $\Theta$ as an operator from $M(\lambda)\otimes
M(\lambda)$ to itself. There is a natural basis for
$M(\lambda)\otimes M(\lambda)$ coming from the PBW basis of
$M(\lambda)$.

\begin{proposition}
The matrix of $\Theta$ acting on $M(\lambda)$ in a PBW basis is $q$--holonomic
with respect  to all the coordinates of the entry and $\lambda$.
\end{proposition}

\begin{proof}
It's enough to prove the statement for each $\Theta_r$.
The result for $\Theta_r$ follows from the fact that the actions
of $e_r^k, f_r^k$ on $M(\lambda)$, as well as $c_k$, are
$q$--holonomic in $k$ and so are all the coordinates of the matrix
entries, by Proposition \ref{Verma} .
\end{proof}

\subsubsection{The $R$--matrix and the braiding} As usual, let us
define the weight  on $M(\lambda)$ by declaring the weight of $F_\n\cdot
\eta$ to be $\lambda - \sum n_i\gamma_i$, where
$\n=(n_1,\dots,n_t)$. The space $M(\lambda)$ is the direct sum of
its weight subspaces.

Let $\DD\co  M(\lambda)\otimes M(\lambda) \to M(\lambda)\otimes
M(\lambda)$ be the linear operator defined by
$$ \DD(x\otimes y) = v^{- (|x|,|y|)} x\otimes y.$$
Clearly $\DD$ is $q$--holonomic; it's called the diagonal part of
the $R$--matrix, which is $R:= \Theta\DD$.

The braiding is $\BB:=R\sigma$, where $\sigma(x\otimes y)=y\otimes
x$. Combining the above results, we get the following:

\begin{theorem}
The entry of the matrix of the braiding acting on $M(\lambda)$ is
$q$--holonomic with respect to all the coordinates and $\lambda$.
\end{theorem}

\begin{remark}
Technically, in order to define the diagonal part $\DD$, one needs to
extend the ground ring to include a $D$-th root of $v$, since
$(\lambda,\mu)$, with $\lambda,\mu \in X$, is in general not an
integer, but belonging to $\frac{1}{D}\Z$.
\end{remark}

\subsubsection{$q$--holonomicity of quantum invariants of knots}
\label{subsub.qqq}

First let us recall the definition of quantum knot invariant.

Using the braiding $\BB\co  M(\lambda)\to M(\lambda)$ one can define
a representation of the braid group $\tau\co  B_m\to
(M(\lambda))^{\otimes m}$ by putting
$$
\tau(\sigma_i):=
\id^{\otimes i-1} \otimes  \BB \otimes \id^{\otimes m-i-1}.
$$
Let $\rho$ denote the half-sum of positive roots. For an element
$x\in \UU$ and an $\UU$--module $V$, the quantum trace is defined
as
$$
\tr_q(x,V) := \tr (xK_{-2\rho}, V).
$$
Suppose a framed knot $\KK$ is obtained by closing a braid
$\beta\in B_m$. We would say that the colored Jones polynomial
is the quantum trace of $\tau(\beta)$. However, since
$M(\lambda)$ is infinite-dimensional, the trace may not make sense.
Instead, we will use a trick of {\em breaking the knot}. Let $\KK'$
denote the {\em long knot} which is the closure of all but the first strand
of $\b$.

Recall that $\tau(\beta)$ acts on $(M(\lambda))^{\otimes m}$. Let
$$
\tau(\beta)(\lambda)_{\n_1,\dots,\n_m}^{\n'_1,\dots,n'_m}\in
\BZ[v^{\pm 1/D}]
$$
be the entries of the matrix
$\tau(\beta)(\lambda)$. We will take partial trace by first
putting $\n_1=\n_1'=0$ and then take the sum over all
$\n_2=\n_2',\dots,\n_m=\n_m'$. The following lemma shows that the
sum is actually finite.

\begin{lemma}
Suppose $\n_1=0$. There are only a finite
number of collections of $(\n_2,\n_3,\dots,\n_m)\in \BN^{t-1}$
such that
$$
\tau(\beta)(\lambda)_{\n_1,\dots,\n_m}^{\n_1,\dots,\n_m}
$$
is not zero.
\end{lemma}

\begin{proof}
Let $M'(\lambda)$ be the maximal proper $\UU$--submodule
of $M(\lambda)$. Then $L(\lambda)= M(\lambda)/M'(\lambda)$ is a
finite dimensional vector space. In particular it has only a
finite number of non-trivial weights. Hence, all except for a
finite number of $\f_\n, \n\in \BN^t,$ are in $M'(\lambda)$.

We present the coefficients $\BB_\pm(\lambda)$ graphically as in
Figure \ref{braiding}.

\begin{figure}[ht!]
$$
\psdraw{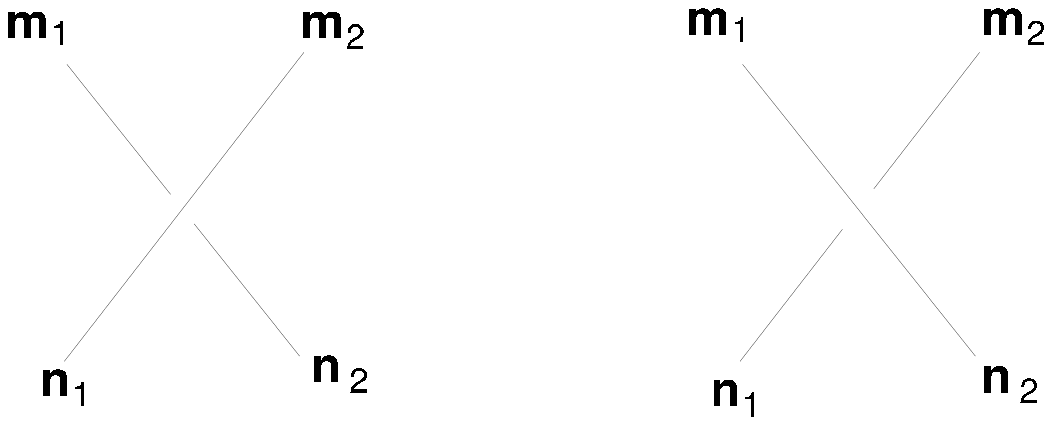}{2.3in}
$$
\caption{$(\BB_+)_{\n_1,\n_2}^{\m_1,\m_2}$ and
$(\BB_-)_{\n_1,\n_2}^{\m_1,\m_2}$}\label{braiding}
\end{figure}

Note that if $(\BB_\pm)_{\n_1,\n_2}^{\m_1,\m_2}$ is not equal to
0, then $\f_{\m_2}$ can be obtained from $\f_{\n_1}$ by action of
an element in $\UU$, and similarly, $\f_{\m_1}$ can be obtained
from $\f_{\n_2}$ by action of an element in $\UU$. Thus if we move
upwards along a string of the braid, the basis element at the top
can always be obtained from the one at the bottom by an action of
$\UU$.

Because the closure of $\beta$ is a knot, by moving around the
braid one can get any point from any particular point. Because
the basis element $\f_0$ is not in $M'(\lambda)$, we conclude
that if
$$\tau(\beta)(\lambda)_{\n_1,\dots,\n_m}^{\n_1,\dots,\n_m}$$
is not 0, with $\n_1=0$, then all the basis vectors $\f_{\n_2},
\dots,\f_{\n_m}$ are not in $M'(\lambda)$, and there are only a
finite number of such collections.
\end{proof}

R
ecall that $2\rho$ is the sum of all positive roots. Let us
define
$$
J_{\KK'}(\lambda)= \sum_{\n_2,\dots\n_m \in \BN^t, \n_1=0,}
(K_{-2\rho}\,\tau(\beta)(\lambda))_{\n_1,\dots,\n_m}^{\n_1,\dots,n_m}.
$$
From $q$--holonomicity of $\tau(\beta)(\lambda)$ it follows that
$J_{\KK'}(\lambda)$ is $q$--holonomic.
$J_{\KK'}(\lambda)$ is a long knot invariant, and is related to the
colored Jones polynomial $J_\KK$ of the knot $\KK$ by
$$
J_\KK(\lambda)= J_{\KK'}(\lambda)\times \dim_q(L(\lambda)),
$$
where $L(\lambda)$ is the finite-dimensional simple $\UU$--module
of highest weight $\lambda$, and $\dim_q(L(\lambda))$ is its
quantum dimension, and is given by the formula
$$ \dim_q(L(\lambda))=
\prod_{\al>0}
\frac{v^{(\lambda+\rho,\al)}-v^{-(\lambda+\rho,\al)}}{v^{(\rho,\al)}
-v^{-(\rho,\al)}}.
$$
Since $\dim_q(L(\lambda))$ is $q$--holonomic in $\lambda$, we see
that $J_\KK(\lambda)$ is $q$--holonomic. This completes the proof
of Theorem \ref{mainn}.
\qed

\begin{remark}
The invariant $J_{\KK'}$ of long knots is sometime more convenient. For
example, $J_\KK(\l)$ might contain fractional power of $q$,
but (if $\KK'$ has framing 0,) $J_{\KK'}(\l)$ is always in
$\BZ[q^{\pm1}]$, see \cite{Le1}. Also the function $J_{\KK'}$ can be
extended to the whole weight lattice.
\end{remark}

\subsection{Proof of Theorem \ref{main}}
\subsubsection{$r_\alpha$ is $q$--holonomic} We will need the
linear maps $r_\alpha, r'_\alpha\co  \UU^{\pm}\to \UU^{\pm}$, as
defined in \cite[Chapter 6]{Jantzen}. Their restriction to
$\UU^-$ is uniquely characterized by the properties:
\begin{equation}
\label{eq.101}
r_\alpha (xy) =  r_\alpha(x)\, y+ v^{(\alpha,|x|)}\, x \,r_\alpha(y)
\qquad
r'_\alpha (xy) =  x \,r'_\alpha(y) + v^{(\alpha,|x|)}\, r'_\alpha(x) \, y
\end{equation}
and for any two fundamental roots $\alpha, \beta$,
(see \cite[Eqn.(6.15.4)]{Jantzen})  and
\begin{equation}
\label{eq.102}
r_\alpha(F_\beta^n) =
r'_\alpha(F_\beta^n) = \d_{\alpha,\beta}\frac{1-v_\al^{2n}}{1-v_\al^2}
F_\alpha^{n-1},
\end{equation}
where $v_\al:= v^{(\al,\al)/2}$; see \cite[Eqn.(8.26.2)]{Jantzen}.

\begin{lemma}
\label{ralpha}
For a fixed $\alpha \in \{\alpha_1,\dots,\alpha_\ell\}$, the
matrix entries of the operators $(r_\alpha)^k, (r'_\alpha)^k :
\UU^-\to \UU^-$ are $q$--holonomic with respect to $k$ and the
coordinates of the matrix entry. Similarly, $(r_\alpha)^k,
(r'_\alpha)^k \co  \UU^+\to \UU^+$ are $q$--holonomic.
\end{lemma}

\begin{proof}
We give a proof for $r_\alpha^k\co  \UU^-\to \UU^-$. The other case
is similar.

There is a
reduced word $w'= s_{i_1}\dots s_{i_t}$ representing the longest
element $\omega_0$ of $W$ such that $\alpha_{i_1}=\alpha$. Then
$w = s_{i_2}\dots s_{i_t} s_{\bar \alpha}$ is another reduced
word representing $\omega_0$, where $\bar \alpha:=-
\omega_0(\alpha)$.

For the  PBW basis of $\UU^-$ associated with $w$ it's known that
$\gamma_t=\alpha$, and thus $f_t= F_\alpha$. According to
\cite[8.26.5]{Jantzen}, for every $x$ in the algebra
$\UU^-[t-1,1]$, one has
$$
r_\alpha(x)=0.
$$
Using Equations (\ref{eq.101}) and (\ref{eq.102}) and induction,
one can easily show that for every $x\in \UU^-[t-1,1]$,
$$
(r_\alpha)^k (f_t^{n_t} x )=
\prod_{i=1}^{k}\frac{1-v^{2n_t-2i+2}_\al}{1-v^2_\al} \, f_t^{n_t-k}x,
$$
This formula, applied to $x= f_{t-1}^{n_{t-1}}\dots f_1^{n_1}$, proves
the statement.
\end{proof}

\subsubsection{$U_q(\msl_2)$ is $q$--holonomic}

\begin{lemma}
Theorem \ref{main} holds true for $\fg=\msl_2$. \label{sl2case}
\end{lemma}

\begin{proof}
The PBW basis for $\UU$ is $F^nK^jE^m$, with $m,n\in
\N$ and $j\in \Z$. First of all we know that
$$
E_\alpha^m F_\alpha^n = \sum_{i=0}^\infty \qbinom{m}{i}_{v_\alpha}
 \qbinom{n}{i}_{v_\alpha}
F_\alpha^{n-i}\, b(K_\alpha;2i-n-m,i)\, E_\alpha^{m-i},
$$
$$
b(K_\alpha;a,i):= \prod_{j=1}^i \frac{K_\alpha v_\alpha^{a-j+1} -
K_\alpha^{-1}v_\alpha^{-a+j-1}}{v_\alpha-v_\alpha^{-1}}.\leqno{\hbox{where}}
$$
Here for any root $\gamma$, one defines  $v_\gamma=
v^{(\gamma,\gamma)/2}$, and $\qbinom{m}{i}_{v_\alpha}$ is the
usual quantum binomial coefficient calculated with $v$ replaced
by $v_\alpha$.

Hence
$$
(F^m K^k E^n )(F^{m'} K^{k'} E^{m'}) = \sum_{i=0}^\infty
F^{m+m'-i} \, a(m,k,n,m',k',n',i)\, E^{n+n'-i},
$$
where
\begin{align*}
a(m,k,n,m',k'&,n',i)\\& = v^{2k(i-m') +
2k'(i-n)}\qbinom{n}{i}\qbinom{m'}{i}[i]! \, b(K;2i-n-m',i)\,
K^{k+k'}.
\end{align*}
The value of the function $a$ is in $\Z[v^{\pm 1}][K^{\pm 1}]$.
Consider the coefficient of $K^r$ in $a$; one gets a function of
$m,n,k,m',n',k',i, r$ with values in $\Z[v^{\pm 1}]$ which is
clearly $q$--holonomic with respect to all variables.  The lemma
follows.
\end{proof}

\subsubsection{$E_\alpha^k, F_\alpha^k\co  \UU \to \UU$ are
$q$--holonomic in $k$}
\label{subsub.EF}

\begin{proposition}
\label{Ealpha}
For a fixed fundamental root $\alpha\in
\{\alpha_1,\dots,\alpha_\ell\}$, the operators $E_\alpha^k,
F_\alpha^k\co  \UU \to \UU$ of left multiplication are $q$--holonomic
with respect to $k$ and all the coordinates of the matrix entry.
Similarly, the right multiplication by $E^k_\alpha, F_\alpha^k$
are $q$--holonomic with respect to $k$ and all the coordinates of the
matrix entry.
\end{proposition}

\begin{proof} (a)\qua Left multiplication by $F_\alpha^k$ and right
multiplication by $E_\alpha^k$.

Choose $w$ as in the proof of Lemma \ref{ralpha}. Then
$f_t=F_\alpha$ and $e_t=E_\alpha$, and an element of the PBW
basis has the form $f_t^{n_t} x K_\beta y e_t^{m_t}$. It's clear
that left multiplication by $F_\alpha$ and right multiplication by
$E_\alpha$  are $q$--holonomic.

(b)\qua Left multiplication by $E_\alpha^k$.

Choose a
reduced word $w= s_{i_1}\dots s_{i_t}$ representing the longest
element $\omega_0$  that begins with $\alpha$:
$\alpha_{i_1}=\alpha$. We have the corresponding PBW basis
$f_i,e_i, i=1,\dots,t$  with $f_1=F_\alpha$ and $e_1=E_\alpha$.
Thus a typical element of the PBW basis has the form
\begin{equation}x F_\alpha^{n_1} K_\beta E_\alpha^{m_1} y,
\label{11}
\end{equation}
where $x= f_t^{n_t} \dots f_2^{n_2}, y=e_2^{m_2}\dots e_t^{m_t}$.
By \cite[8.26.6]{Jantzen}, since  $x\in \UU^-[t,2]$, one has
$r'_\alpha(x)=0$. Using formula \cite[6.17.1]{Jantzen}, one can
easily prove by induction that
$$ (E_\alpha)^k x = \sum_{i=0}^\infty v^{i-ik}\qbinom{k}{i}_{v_\alpha}
\frac{K_\alpha^i}{(v_\alpha-v_\alpha^{-1})^i}(r_\alpha)^i(x)
E_\alpha^{k-i}.$$
Using this formula one can move the $E_\alpha$ past $x$ in the
expression (\ref{11}), (there appear $r_\alpha$ and $K_\alpha$),
then one moves $E_\alpha$ past $F_\alpha$ using the $\msl_2$ case.
The last step is moving past $K_\beta$ is easy, since
$$E_\alpha K_\beta = v^{-(\beta,\alpha)} K_\beta E_\alpha.$$
 Using Lemmas \ref{ralpha} and \ref{sl2case}, we see that each ``moving step"
 is $q$--holonomic.  Hence we get the
result for the left multiplication by $E_\alpha^k$.

(c)\qua Right multiplication by $F_\alpha^k$.

The proof is similar. We use the same basis (\ref{11}) as
in the case b).  For $y$, by Lemma 8.26 of \cite{Jantzen}, one
has $r_\alpha(y)=0$. Hence using induction based on the formula
(6.17.2) of \cite{Jantzen} one can show that
$$ y F_\alpha^n = \sum_{i=0}^\infty
\frac{v_\alpha^{i(n-i)}}{(v_\alpha^{-1}-v_\alpha)^i}\qbinom{n}{i}_{v_\alpha}
F_\alpha^{n-i} K_\alpha^{-i} (r'_\alpha)^i(y).$$
Using this formula, and the results for $r'_\alpha$ (Lemma
\ref{ralpha}) and $\msl_2$ (Lemma \ref{sl2case}) we can move
$F_\alpha$ to the right.
\end{proof}

\subsubsection{$T_\alpha$ is $q$--holonomic}
\label{subsub.talpha}

\begin{proposition}
For a fixed fundamental root $\alpha\in
\{\alpha_1,\dots,\alpha_\ell\}$, the braid operator $T_\alpha\co 
\UU \to \UU$ and its inverse $T_\alpha^{-1}$ are $q$--holonomic.
\label{Talpha}
\end{proposition}

\begin{proof}
By Proposition \ref{prop.trans} we can use any PBW basis.

Choose a reduced word $w'= s_{i_1}\dots s_{i_t}$ representing the
longest element $\omega_0$  that begins with $\alpha$:
$\alpha_{i_1}=\alpha$.  Then $w= s_{i_2}\dots s_{i_t} s_{\bar
\alpha}$ is another reduced word representing $\omega_0$, where
$\bar\alpha$ is the dual of $\alpha$: $\bar \alpha=
-\omega_0(\alpha)$.

We use $f_r$ to denote $f_r(w)$, and $f'_r$ to denote $f_r(w')$.
The relation between the two PBW basis of $w$ and $w$ is as
follows: For $1\le r\le t-1$,
$$ T_\alpha(f_r) = f'_{r+1}, \qquad T_\alpha(e_r) =
e'_{r+1}.$$
Besides, $f_t=F_\alpha = f_1', e_t= E_\alpha= e_1'$.

We will consider the matrix entry of $T_\alpha\co  \UU\to \UU$ where
the source space is equipped with the PBW corresponding to $w$, while the
target space with the PBW basis corresponding to $w'$.

From \cite[Chapter 8]{Jantzen}, recall that:
$$
T_\alpha(F_\alpha) = -K_\alpha^{-1} E_\alpha, \qquad
T_\alpha(E_\alpha) = -F_\alpha
K_\alpha.
$$
Hence
$$
T_\alpha(F_\alpha^n) = (-1)^n v_\alpha^{n(n-1)} K_\alpha^{-n}
E_\alpha^n, \qquad T_\alpha(E_\alpha^m) = (-1)^m
v_\alpha^{-m(m-1)} F_\alpha^{m} K_\alpha^m.
$$
For a basis element $x_\s=f_t^{n_t}\dots f_1^{n_1} K_\beta
e_1^{m_1}\dots e_t^{m_t}$,  we have
\begin{align*}
T_\alpha(x_\s) = d_\alpha(n_t,m_t)&\,
 K_\alpha^{-n_t}
E_\alpha^{n_t} \times (f_t')^{m_{t-1}}\\& \dots (f_1')^{n_2}
K_{s_\alpha\beta} (e_1')^{m_2}
 \dots (e'_t)^{m_{t-1}}\times
F_\alpha^{m_t} K_\alpha^{m_t},
\end{align*}
$$
d_\alpha(n_t,m_t) :=
(-1)^{n_t+m_t}v_\alpha^{n_t(n_t-1)-m_t(m_1-1)}.\leqno{\hbox{where}}
$$
The left or right multiplication by $K_\alpha^n$ is
$q$--holonomic with respect to $n$ and all the coordinates. The left
multiplication by $E^{n_t}$, as well as the right multiplication
my $F_\alpha^{m_t}$ is $q$--holonomic with respect to $n_t$ and all
coordinates, by Proposition \ref{Ealpha}. One then can conclude
that $T_\alpha$ is $q$--holonomic.

The proof for $T_\alpha^{-1}\co \UU\to \UU$ is similar. One should
use the PBW basis of $w'$ for the source, and that of $w$ for the
target.
\end{proof}

\subsubsection{Proof of Theorem \ref{main}}
\label{proof}

It is clear that for each $\j\in \Z^\ell$, the operator $\K^\j\co  \UU
\to \UU$ of left multiplication is $q$--holonomic.

Fix a reduced word $w$ representing the longest element of $W$.
It suffices to show that for each $1\le r\le t$ the operators
$e_r^k, f_r^k\co  \UU\to \UU$ (left multiplication) are $q$--holonomic
with respect to all variables, including $k$.

This is true if $e_r=E_\alpha$ and $f_r=F_\alpha$, where $\alpha$
is one of the fundamental roots, by Proposition \ref{Ealpha}. But
any $e_r$ or $f_r$ can be obtained from $E_\alpha$ and $F_\alpha$
by actions of product of various $T_{\alpha_i}$'s. Hence from
Proposition \ref{Talpha} we get Theorem \ref{main}.
\qed

\appendix

\section{Appendix: Proof of Proposition \ref{prop.trans} for $A_2$ and $B_2$}
\label{sec.appendix}

In this appendix we will prove Proposition \ref{prop.trans}
for the rank $2$ Lie algebras $A_2$ and $B_2$. We will achieve this by a
brute-force calculation.

First, let us discuss some simplification, due to symmetry. The
transition matrix of $\UU$ leaves invariant each of
$\UU^+,\UU^-,\UU^0$. On $\UU^0$ the transition matrix is
identity. Hence it's enough to consider the restriction of the
transition matrix in $\UU^-$ and $\UU^+$. Furthermore, the Cartan
symmetry (the operator $\tau$ of \cite{Jantzen}) reduces the case
of $\UU^+$  to that of $\UU^-$.

\subsection{The case of $A_2$}
\label{sub.A2}

There are two fundamental roots denoted by $\alpha$ and $\beta$. The set
of positive roots is $\{\a,\b,\a+\b\}$.
The reduced representations of the
longest element of the Weyl group are $w=s_1 s_2 s_1$ and $w'=s_2
s_1s_2$, where $s_1=s_{\a}$ and $s_2=s_{\b}$.

The total ordering (see Section \ref{subsub.ordering}) of the set of positive
roots corresponding to $w$ and $w'$ are, respectively:
\begin{eqnarray*}
(\ga_1,\ga_2,\ga_3) &=& (\a,\a+\b,\b) \\
(\ga_{1'},\ga_{2'},\ga_{3'}) &=& (\b,\a+\b,\a).
\end{eqnarray*}
Notice that $\ga_{i'}=\ga_{3-i}$ for $i=1,\dots,3$.

The PBW basis of $\UU^-$ (see Section \ref{subsub.PBW})
corresponding to $w$ and $w'$ are, respectively:
$$
\{f_3^m f_2^n f_1^p \, | \,\,  m,n,p\in \BN \}, \qquad
\{f_{3'}^m f_{2'}^n f_{1'}^p \, | \,\,  m,n,p\in \BN \},
$$
where
\begin{eqnarray*}
(f_3,f_2,f_1) &=&  (F_\beta, T_\alpha(F_\beta)= -v [F_\beta,F_\alpha ]_q =
F_\beta F_\alpha -v F_\alpha F_\beta, F_\alpha) \\
(f_{3'},f_{2'},f_{1'}) &=& (F_\alpha , T_\beta(F_\alpha)=
F_\alpha F_\beta -v F_\beta F_\alpha,
F_\beta ).
\end{eqnarray*}
From explicit formulas of \cite[section 5]{Lusztig2} it follows that:

\begin{lemma}
\label{1}
The structure constants of $\UU^-$, in the basis of $w$, is $q$--holonomic.
\end{lemma}

Let us define a scalar product $(\cdot, \cdot)$
on $\UU^-$ such that the PBW basis of $w$ is an orthonormal basis. Since
$$
f_{3'}^{m'} f_{2'}^{n'} f_{1'}^{p'}
=\sum_{m,n,p}
(f_{3'}^{m'} f_{2'}^{n'} f_{1'}^{p'}, f_3^m f_2^n f_1^n)
f_3^m f_2^n f_1^p
$$
Proposition \ref{prop.trans} is equivalent to showing that
$$
(f_{3'}^{m'} f_{2'}^{n'} f_{1'}^{p'}, f_3^m f_2^n f_1^n)
$$
is $q$--holonomic in all variables $m,n,p,m',n',p'$.

Since multiplication is $q$--holonomic in the  PBW basis of $w$ (see
Lemma \ref{1}), it suffices to show that
$$(f_{i'}^{k}, f_3^m f_2^n f_1^n)$$
is $q$--holonomic in $k,m,n,p$ for each $i=1,2,3$. This is clear for
$i=1$ or $i=3$, since $f_{1'}=f_3$ and $f_{3'}=f_1$. As for $f_{2'}$,
an easy induction shows that
$$
f_{2'}^n = (-v)^{-n} \sum_{k=0}^\infty v^{-k(k-3)/2}
(v-v^{-1})^k \qbinom{n}{k} f_3^k f_2^{n-k} f_1^k.
$$
and the statement also holds true for $i=2$. This proves Proposition
\ref{prop.trans} for $A_2$.

\subsection{The case of $B_2$}
\label{sub.B2}

There are two fundamental roots denoted here by $\alpha$ and
$\beta$, where $\alpha$ is the short root. The set of positive
roots is $\{\a,\b,2\a+\b,\a+\b\}$. The reduced representations
of the longest element of the Weyl group are $w=s_1 s_2 s_1 s_2$
and $w'=s_2 s_1s_2 s_1$, where $s_1=s_{\a}$ and $s_2=s_{\b}$.

The total ordering of the set of positive
roots corresponding to $w$ and $w'$ are, respectively:
\begin{eqnarray*}
(\ga_1,\ga_2,\ga_3,\ga_4) &=& (\a,2 \a+\b,\a+\b,\b) \\
(\ga_{1'},\ga_{2'},\ga_{3'},\ga_{4'}) &=& (\b,\a+\b,2\a+\b,\a).
\end{eqnarray*}
Notice that $\ga_{i'}=\ga_{4-i}$ for $i=1,\dots,4$.

The PBW basis of $\UU^-$ (see Section \ref{subsub.PBW})
corresponding to $w$ and $w'$ are, respectively:
$$
\{f_4^l f_3^m f_2^n f_1^p \, | \,\,  l,m,n,p\in \BN \},
\qquad
\{f_{4'}^{l} f_{3'}^m f_{2'}^n f_{1'}^p \, | \,\,  l,m,n,p\in \BN \},
$$
where
\begin{eqnarray*}
(f_4,f_3,f_2,f_1) &=& (F_\beta ,
F_\beta F_\alpha-v^2 F_\alpha F_\beta ,
\frac{F_\beta F_\alpha^2}{[2]} -v F_\alpha F_\beta
F_\alpha + \frac{v^2 F_\alpha^2 F_\beta}{[2]}, F_\alpha)  \\
(f_{4'}, f_{3'},f_{2'},f_{1'}) &=& (F_\alpha ,
\frac{v^2 F_\beta F_\alpha^2}{[2]} -v F_\alpha F_\beta F_\alpha +
\frac{F_\alpha^2 F_\beta}{[2]},
F_\alpha F_\beta -v^2 F_\beta F_\alpha, F_\beta ).
\end{eqnarray*}
It follows from \cite{Lusztig2} that:

\begin{lemma}
\label{2}
The structure constants of $\UU^-$, in the basis of $w$, is $q$--holonomic.
\end{lemma}

Let us define a scalar product $(\cdot,\cdot)$ on $\UU^-$ such that the
PBW basis of $w$ is an orthonormal basis. Then Proposition \ref{prop.trans}
is equivalent to
$$(f_{4'}^{l'} f_{3'}^{m'} f_{2'}^{n'} f_{1'}^{p'}, f_4^lf_3^m f_2^n f_1^n)$$
is $q$--holonomic in all variables $l,m,n,p,l'm',n',p'$.

Since multiplication is $q$--holonomic in the  PBW basis of $w$ (see
Lemma \ref{2}), it suffices to show that
$$(f_{i'}^{k}, f_4^l f_3^m f_2^n f_1^p)$$
is $q$--holonomic in $k,l,m,n,p$ for each $i=1,2,3,4$. This is
clear for $i'=1$ or $i'=4$, since $f_{1'}=f_4$ and $f_{4'}=f_1$.
As for $i'=2$ and $i'=3$, the formula of \cite[Section 37.1]{Lusztig}
shows that
\begin{eqnarray*}
f_{2'}^n &=& \sum_{i=0}^{n} (-1)^i \frac{v^{2i} F_\beta^{i} F_\alpha^n
F_\beta^{n-i} } {[n-i]_\beta![i]_\beta!} \\
f_{3'}^n &=& \sum_{i=0}^{2n} (-1)^i
\frac{v^i F_\alpha^{2n-i} F_\beta^n F_\alpha^i}{[2n-i]![i]!}
\end{eqnarray*}
and since $F_{\a}=f_{4'}$ and $F_{\b}=f_{1'}$,
the cases of  $i'=2'$ and $i'=3'$ reduce to the cases of
$i'=1'$ and $i'=4'$. This proves Proposition \ref{prop.trans} for $B_2$.

\begin{remark}
\label{rem.appG2}
If Lemma \ref{1} holds for $G_2$, then we can prove Proposition
\ref{prop.trans} for $G_2$.
\end{remark}


\begin{thebibliography}{[EMSS]}


\bibitem{BN} \textbf{D Bar-Natan},
        {\em ColoredJones.nb}, Mathematica program, part of ``KnotAtlas",
        April 2003.

\bibitem{B1}
\textbf{I\,N Bern{\v{s}}te{\u\i}n}, \emph{Modules over a ring of differential
  operators. {A}n investigation of the fundamental solutions of equations with
  constant coefficients}, Funkcional. Anal. i Prilo\v zen. 5 (1971) 1--16,
  English translation: 89--101
  \MR{0290097}

\bibitem{B2}
\textbf{I\,N Bern{\v{s}}te{\u\i}n}, \emph{Analytic continuation of
  generalized functions with respect to a parameter}, Funkcional. Anal. i
  Prilo\v zen. 6 (1972) 26--40, 
  English translation: 273--285 \MR{0320735}

\bibitem{Bj}
\textbf{J-E Bj{\"o}rk}, \emph{Rings of differential operators}, volume~21 of
  \emph{North-Holland Mathematical Library}, North-Holland Publishing Co.,
  Amsterdam (1979)\relax \MR{549189}

\bibitem{Bo}
\textbf{A Borel}, \textbf{P-P Grivel}, \textbf{B Kaup}, \textbf{A Haefliger},
  \textbf{B Malgrange}, \textbf{F Ehlers}, \emph{Algebraic {$D$}-modules},
  Perspectives in Mathematics 2, Academic Press, Boston,
  MA (1987)\relax \MR{882000}

\bibitem{C}
\textbf{P Cartier}, \emph{D\'emonstration ``automatique'' d'identit\'es et
  fonctions hyperg\'eom\'etriques (d'apr\`es {D}.\ {Z}eilberger)}, Ast\'erisque
   (1992) Exp.\ No.\ 746, 3, 41--91\relax \MR{1206064}

\bibitem{Cou}
\textbf{S\,C Coutinho}, \emph{A primer of algebraic {$D$}-modules}, volume~33
  of \emph{London Mathematical Society Student Texts}, Cambridge University
  Press, Cambridge (1995)\relax \MR{1356713}

\bibitem{CS}
\textbf{F Chyzak}, \textbf{B Salvy}, \emph{Non-commutative
  elimination in {O}re algebras proves multivariate identities}, J. Symbolic
  Comput. 26 (1998) 187--227\relax \MR{1635242}

\bibitem{GLZ} \textbf{S Garoufalidis}, \textbf{T\,T\,Q Le}, \textbf{D Zeilberger},
        {\em The quantum MacMahon Master Theorem},
        \arxiv{math.QA/0303319}, to appear in Proc. Nat. Acad. Sci.

\bibitem{G1} \textbf{S Garoufalidis},
        {\em Difference and differential equations for the colored Jones
        function}, \arxiv{math.GT/0306229}

\bibitem{G2} \textbf{S Garoufalidis},
        {\em On the characteristic and deformation varieties of a knot},
        from: ``Proceedings of the CassonFest (Arkansas and Texas 2003)",
        \gtmref7{2004}{12}{291}{309}

\bibitem{GaS} \textbf{S Garoufalidis}, \textbf{X Sun},
        {\em The $C$--polynomial of a knot}, preprint (2005)
        \arxiv{math.GT/0504305}

\bibitem{Gelca1}
\textbf{R Gelca}, \emph{Non-commutative trigonometry and the
  {$A$}-polynomial of the trefoil knot}, Math. Proc. Cambridge Philos. Soc. 133
  (2002) 311--323\relax \MR{1912404}

\bibitem{Gelca2}
\textbf{R Gelca}, \textbf{J Sain}, \emph{The noncommutative
  {A}-ideal of a {$(2,2p+1)$}-torus knot determines its {J}ones polynomial}, J.
  Knot Theory Ramifications 12 (2003) 187--201\relax \MR{1967240}

\bibitem{Habiro}
\textbf{K Habiro}, \emph{On the quantum {$\rm sl\sb 2$} invariants of knots
  and integral homology spheres}, from: ``Invariants of knots and 3-manifolds
  (Kyoto, 2001)'', \gtmref4{2002}5{55}{68} \MR{2002603}

\bibitem{HL} \textbf{K Habiro}, \textbf{T\,T\,Q Le},
        in preparation

\bibitem{Jantzen}
\textbf{J Jantzen}, \emph{Lectures on quantum groups}, volume~6 of
  \emph{Graduate Studies in Mathematics}, American Mathematical Society,
  Providence, RI (1996)\relax \MR{1359532}


\bibitem{Soibelman}
\textbf{L\,I Korogodski}, \textbf{Y\,S Soibelman}, \emph{Algebras of
  functions on quantum groups. {P}art {I}}, volume~56 of \emph{Mathematical
  Surveys and Monographs}, American Mathematical Society, Providence, RI
  (1998)\relax \MR{1614943}

\bibitem{Hi}
\textbf{K Hikami}, \emph{Difference equation of the colored {J}ones
  polynomial for torus knot}, Internat. J. Math. 15 (2004) 959--965\relax
  \MR{2106155}

\bibitem{Le1}
\textbf{T\,T\,Q Le}, \emph{Integrality and symmetry of quantum link
  invariants}, Duke Math. J. 102 (2000) 273--306\relax \MR{1749439}

\bibitem{Le2} \textbf{T\,T\,Q Le}, 
        {\em The Colored Jones Polynomial and the $A$--Polynomial of 
        Knots}, \arxiv{math.GT/0407521}

\bibitem{Lusztig}
\textbf{G Lusztig}, \emph{Introduction to quantum groups}, volume 110 of
  \emph{Progress in Mathematics}, Birkh\"auser Boston Inc., Boston, MA
  (1993)\relax \MR{1227098}

\bibitem{Lusztig2}
\textbf{G Lusztig}, \emph{Quantum groups at roots of {$1$}}, Geom.
  Dedicata 35 (1990) 89--113\relax \MR{1066560}

\bibitem{Ml}
\textbf{B Malgrange}, \emph{\'{E}quations diff\'erentielles \`a coefficients
  polynomiaux}, volume~96 of \emph{Progress in Mathematics}, Birkh\"auser
  Boston Inc., Boston, MA (1991)\relax \MR{1117227}

\bibitem{Ma}
\textbf{G Masbaum}, \emph{Skein-theoretical derivation of some formulas of
  {H}abiro}, \agtref3{2003}{17}{537}{556}
  \MR{1997328}

\bibitem{PR1}
\textbf{P Paule}, \textbf{A Riese}, \emph{A {M}athematica {$q$}-analogue
  of {Z}eilberger's algorithm based on an algebraically motivated approach to
  {$q$}-hypergeometric telescoping}, from: ``Special functions, $q$--series and
  related topics (Toronto, ON, 1995)'', Fields Inst. Commun. 14, Amer. Math.
  Soc., Providence, RI (1997)  179--210\relax \MR{1448687}

\bibitem{PR2} \textbf{P Paule}, \textbf{A Riese}, Mathematica software,
 available at:\newline {\footnotesize{\url{http://www.risc.uni-linz.ac.at/research/combinat/risc/software/qZeil/}}}

\bibitem{PWZ}
\textbf{M Petkov{\v{s}}ek}, \textbf{H\,S Wilf}, \textbf{Doron
  Zeilberger}, \emph{{$A=B$}}, A K Peters Ltd., Wellesley, MA (1996)\relax
  \MR{1379802}

\bibitem{RT}
\textbf{N\,Yu Reshetikhin}, \textbf{V\,G Turaev}, \emph{Ribbon graphs and their
  invariants derived from quantum groups}, Comm. Math. Phys. 127 (1990)
  1--26\relax \MR{1036112}

\bibitem{R} \textbf{A Riese},
        {\em qMultisum--A package for proving $q$--hypergeometric multiple
        summation identities},
        preprint (2002)

\bibitem{S}
\textbf{C Sabbah}, \emph{Syst\`emes holonomes d'\'equations aux
  {$q$}-diff\'erences}, from: ``$D$--modules and microlocal geometry (Lisbon,
  1990)'', de Gruyter, Berlin (1993)  125--147\relax \MR{1206016}

\bibitem{Tu}
\textbf{V\,G Turaev}, \emph{The {Y}ang-{B}axter equation and invariants of
  links}, Invent. Math. 92 (1988) 527--553\relax \MR{939474}

\bibitem{WZ}
\textbf{H\,S Wilf}, \textbf{D Zeilberger}, \emph{An algorithmic proof
  theory for hypergeometric (ordinary and ``{$q$}'') multisum/integral
  identities}, Invent. Math. 108 (1992) 575--633\relax \MR{1163239}

\bibitem{Y}
\textbf{L Yen}, \emph{A two-line algorithm for proving {$q$}-hypergeometric
  identities}, J. Math. Anal. Appl. 213 (1997) 1--14\relax \MR{1469359}

\bibitem{Z}
\textbf{D Zeilberger}, \emph{A holonomic systems approach to special
  functions identities}, J. Comput. Appl. Math. 32 (1990) 321--368\relax
  \MR{1090884}




\end{thebibliography}
\end{document}